\documentclass[a4papger,10pt]{article}
\usepackage{booktabs}%表の加工
\usepackage{lscape}
\usepackage{latexsym,indentfirst,ascmac}
\usepackage{hhline}
\usepackage{ascmac}%見出しつきの枠  
\usepackage{bm,amsmath,amssymb,amsthm,latexsym,indentfirst,amsfonts}
\usepackage{enumerate}  
\usepackage{fancyhdr}
\usepackage{verbatim}    
\usepackage{makeidx}
\usepackage{here}
\usepackage{subcaption}
\usepackage{caption}
\usepackage{mathrsfs} 
\usepackage{multirow}
\usepackage[dvipdfmx]{graphicx,color}
\numberwithin{equation}{section}

\renewcommand{\theenumi}{(A-\arabic{enumi})}

\newtheorem{theorem}{Theorem}

\newtheorem{lemma}{Lemma}
 
\newtheorem{lemma*}{Lemma A.}
\newtheorem{remark}{Remark} 

\newtheorem{definition}{Definition} 
\newtheorem{assumption}{Assumption}

\newtheorem{corollary}{Corollary}

\newsavebox{\circlebox}
\savebox{\circlebox}{\fontencoding{OMS}\selectfont\scriptsize\char13}%\Large\char13}
\newlength{\circleboxwdht}
\title{Time series quantile regression using random forests}

\author{%\color{red}
Hiroshi Shiraishi\footnote{Department of Mathematics, Keio University,%, Yokohama, Kanagawa, Japan. 
E-mail: shiraishi@math.keio.ac.jp.},\quad
Tomoshige Nakamura\footnote{Department of Mathematics, Keio University},\quad
Ryotato Shibuki\footnote{Department of Mathematics, Keio University}
}

\date{}
\begin{document} 
%\gtfamily
\maketitle
%%%%%%    TEXT START    %%%%%%
\begin{abstract}
We discuss an application of Generalized Random Forests (GRF) proposed by Athey et al. (2019)
%\cite{athey2019generalized}
 to quantile regression for time series data. 
We extracted the theoretical results of the GRF consistency for i.i.d. data to time series data. 
In particular, in the main theorem, based only on the general assumptions for time series data in Davis and Nielsen (2020)
%\cite{davis2020modeling}
, and trees in Athey et al. (2019)
%\cite{athey2019generalized}
, we show that the tsQRF (time series Quantile Regression Forests) estimator is consistent. 
Davis and Nielsen (2020) 
%\cite{davis2020modeling}
 also discussed the estimation problem using Random Forests (RF) for time series data, but the construction procedure of the RF treated by the GRF is essentially different, and different ideas are used throughout the theoretical proof. 
In addition, a simulation and real data analysis were conducted. 
In the simulation, the accuracy of the conditional quantile estimation was evaluated under time series models. 
In the real data using the Nikkei Stock Average, our estimator is demonstrated to be more sensitive than the others in terms of volatility, thus preventing underestimation of risk. 
\end{abstract}

\section{Introduction}\label{sec1}
    Quantile Regression (QR), proposed by Koenker and Bassett (1978)
    %\cite{koenker1978regression}
    , is regression model that has been applied in various fields. 
QR analyzes the effects of covariates on outcomes by focusing on quantiles rather than means. 
    Therefore, it can flexibly analyze the effect of covariates on the tail of the conditional distribution, which cannot be captured by regression on the mean. QR is used in a wide range of fields, including economics, medicine, and epidemiology, and is applicable not only to cross-sectional data, but also to panel data, for which the theory is well developed (Koenker , 2005)
    %\cite{koenker2005})
    . 
    In the analysis of time series data, many studies have focused on the dynamics of the mean of the series, and on that of the conditional distribution. 
    Therefore, the estimation of conditional quantiles of time series data by QR can capture the local and quantile-specific dynamics of time series and is expected to improve the quality of analysis.
    
    Research on QR for time series data is ongoing. 
Value-at-Risk (VaR) is one of the criteria used to measure the market risk of an asset in the field of risk management. 
    Because VaR is defined as the quantile of the asset return at time $t$, it is necessary to estimate the conditional quantile regression function to estimate VaR. 
    There are two methods to estimate the quantile: parametric and nonparametric methods. 
    For parametric QR, Koul and Saleh (1995) 
    %\cite{Koul1995}
     proposed a conditional quantile estimator for AR models, Koenker and Zhao (1996) 
    %\cite{koenker_and_zhao_1996}
     proposed  ARCH models, Taylor (1999) 
    %\cite{taylor1999quantile}
     introduced a linear VaR model, and Chernozhukov and Umantsev (1996)
    %\cite{chernozhukov2001conditional}
     introduced a quadratic VaR model. 
   The CAViaR model proposed by Engle and Manganelli (2004) 
   %\cite{engle2004caviar}
    is a broader class of models in time series QR for VaR estimation. 
    In CAViaR, the effect of past VaR on VaR at time $t$ is described as linear, and the effect of the observed time series on VaR at time $t$ is modeled linearly or nonlinearly by the  researcher. 
    Parametric QR has good properties in terms of interpretation and ease of implementation.
    However, parametric models have a serious bias if they are misspecified. 
    To avoid this problem, Hall (1999)
    %\cite{hall1999methods}
    , Cai (2002)
    %\cite{cai2002regression}
    , Wu et al. (2007) 
    %\cite{Wu_Yu_Mitra2007}
     and Cai and Wang (2008)
    %\cite{Cai2008}
     proposed a nonparametric method for estimating the conditional distribution function using kernel smoothing. 
     In particular, Cai (2002)
     %\cite{cai2002regression}
      showed that the Weighted Nadaraya-Watson (WNW) estimator proposed by Hall (1999) 
     %\cite{hall1999methods}
      for time series satisfying the $\alpha$-mixing property has consistency and asymptotic normality. 
     Furthermore, Cai and Wang (2008) 
     %\cite{Cai2008}
      proposed a Weighted Double-Krnel Local Linear (WDKLL) estimator, an extension of the WNW estimator, and showed its consistency and asymptotic normality. 
     However, in a nonparametric quantile regression estimator using smoothing with kernel functions, the accuracy of the estimator highly depends on the choice of kernel function and  bandwidth parameter.  
    
    Random forest is a representative algorithm in machine learning that has been successfully used in various applications since its proposal in Breiman (2001)
    %\cite{breiman2001random}. 
    In recent years, the asymptotic properties of the estimators obtained by random forests have been studied in terms of consistency (Biau et al. (2008)
    %\cite{biau2008consistency}
    ; Denil et al. (2013)
    %\cite{denil2013consistency}
    ; Scornet et al. (2015)
    %\cite{scornet2015consistency})
    , and asymptotic normality (Wager and Athey (2018)
    %\cite{wager2018estimation}
    ). 
    Consequently, random forests are now treated not only as a predictive model but also as a nonparametric statistical model.
    
    Various extensions of random forests have also been proposed. 
    A method for estimating conditional quantiles using random forests for i.i.d. data is quantile regression forests (Meinshausen (2006)
    %\cite{meinshausen2006quantile}
    ). 
    Davis and Nielsen (2020) 
    %\cite{davis2020modeling}
     showed that random forest estimators are consistent with the problem of regression on the conditional mean of time series data with $\alpha$-mixing properties.
    
    Among the recent extensions of random forests, the most notable is the Generalized Random Forests (GRF) by Athey et al. (2019)
 %\cite{athey2019generalized}
 , which estimates a parameter defined as the solution to a local estimating equation. 
    Athey et al. (2019) 
    %\cite{athey2019generalized}
     showed that the estimator obtained by the GRF had consistency and asymptotic normality under i.i.d. data observations. 
    Thus, from a theoretical perspective, the GRF can be used to estimate parameter. 
    
    In this study, we devote quantile regression estimators to time series data using GRF.
    Thus, it is necessary to extend the existing studies to the following points: 
    Davis and Nielsen's results for time series data focus on estimating a regression model for the mean, which is applicable to quantile regression. 
    However,  this method assumes that the number of samples in each terminal node (leaf) is proportional to sample size. 
    It is difficult to verify this assumption when considering applications. 
The subsample size included in the leaf should not depend on the sample size. 
   In the GRF, the subsample size in the leaf can be fixed to resolve the above problem.  
   The contribution of this study is to propose a time series Quantile Regression Forest (tsQRF) using the GRF framework for $\alpha$-mixing time series data. 
    Furthermore, we extend the theoretical results of GRF consistency for i.i.d. data to time series data. 
    In particular, in the main theorem, based only on the general assumptions for time series data in Davis and Nielsen (2020)
    %\cite{davis2020modeling}
     and trees in Athey et al. (2019) 
    %\cite{athey2019generalized}
    , we showed that the tsQRF estimator is consistent.
    We also visualized the convergence of the tsGRF estimator through several simulation settings and compared its conditional quantile estimation accuracy with that of the WNW estimator. 
    Furthermore, we fit the proposed method to the Nikkei Stock Average data and compared it with the WNW estimator to clarify the high sensitivity of tsQRF to time series volatility. 
    
    The remainder of this paper is organized as follows. 
    Section 2 discusses the properties of the quantile regression estimator using the GRF. 
    Therefore, we define the double-sample tree score and GRF score, and then show the consistency of our estimator.    
    Section 3 discusses the asymptotic behavior of the estimator and compares its accuracy with that of the WNW estimator through simulations.  
    Section 4 illustrates the result of applying the proposed method to Nikkei Stock Average data and compares the results with those obtained using the WNW estimator.  
    Finally, Section 5 summarizes the discussion and presents future issues. 
    Concrete proofs of the theoretical results presented in Section 2 are provided in the appendix.

\section{Theoretical Results}\label{sec2}
Let $\{\varepsilon_t\}_{t\ge1}$ be a sequence of i.i.d. random variables with $\mathbb E[\varepsilon_t]=0$ and $\mathbb{E}[\varepsilon_t^2]<\infty$, and fix an integer $p\ge1$. 
Given a measurable function $g:\mathbb R^p\to \mathbb R$, define a process $\{Y_t\}_{t\ge 1}$
\begin{align}
    Y_t = g(X_t)+\varepsilon_t,\quad X_t=(Y_{t-1},\ldots,Y_{t-p}),\quad t\ge1. \label{eq2-1}
\end{align}
In this paper, we impose the followings. 
\begin{assumption}
\begin{enumerate}
    \item \label{A1} The random variable $\varepsilon_1$ admits a density $f_{\varepsilon}$ which is positive almost everywhere on $\mathbb R$ and, for some $c\in (0,\infty)$,
	\[ \mathbb{E}[|\varepsilon_1|^m] \le m!c^{m-2},\quad m=3,4,\ldots \]
	Moreover, the cumulative distribution function $F_{\varepsilon}$ of $\varepsilon_1$ satisfies 
	\[ \sup_{x\in \mathbb R} \frac{F_{\varepsilon}(x+\lambda) }{F_{\varepsilon}(x)} <\infty \]   
	for any $\lambda\in (0,\infty)$.  
	\item \label{A2} The function $g$ in (\ref{eq2-1}) is bounded and Lipschitz continuous.
\end{enumerate}
\end{assumption}
\begin{remark}
%The %Assumptions
 \ref{A1} and \ref{A2} include assumptions (A1) and (A2) of Davis and Nielsen (2020) %\cite{davis2020modeling}
, which implies that the process $\{Y_t\}_{t\ge 1}$ is strictly stationary and the $p$-th order Markov chain strictly ensured geometrical ergodicity (c.f., An and Huang (1996) %\cite{an1996geometrical}
,Theorem3.1) and exponentially $\alpha$-mixing (c.f., Doukhan (2012) %\cite{doukhan2012mixing}
,p.89). 
Hence, let $\alpha(n)$ denote the $\alpha$-mixing coefficient. 
Then we have $\alpha(n)\lesssim e^{-n}$\footnote{For two sequence $\{a_n\}_{n\ge 1}$ and $\{b_n\}_{n\ge 1}$, we write $a_n\lesssim b_n$ if there exists a constant $c\ge 1$ such that $a_n\le c b_n$ for all $n$.} under %assumptions
 \ref{A1} and \ref{A2}. 
\end{remark}

\begin{remark}
On the other hand, we do not assume assumption (A-3) of Davis and Nielsen (2020) %\cite{davis2020modeling}
 which is a condition for the minimum subsample size falling in each leaf. 
Instead, we introduce PNN (Potentical Nearest Neighber) $k$-set in the splitting rule following Wager and Athey (2018) %\cite{wager2018estimation}
 and Athey et.al. (2019) %\cite{athey2019generalized}
  (see %Assumption 
\ref{A5}%(A-5)
). 
\end{remark}

Let $\mathcal Y\subset \mathbb R$ and $\mathcal X\subset \mathbb R^p$ be compact subsets of the spaces taken value of $Y_t$  and $X_t$. 
Let $\mathcal Q\subset \ell^{\infty}(\mathcal X)$\footnote{The $\ell^{\infty}(\mathcal X)$ is a set of all uniformly bounded real functions on $\mathcal X$.} be a set of function 
$$\mathcal Q :=\{q: \mathcal X \to \mathcal Y\}.$$
In this paper, under some fixed $\tau \in (0,1)$, we interested in the estimation or prediction of  conditional $\tau$-quantile function $q_0:=(q_0(x))_{x\in \mathcal X}\in \mathcal Q$ defined as solution of the following locall estimationg equation: 
\begin{align}
\Psi(q)(x) := \mathbb E[\psi_{q_0(x)}(Y_t)| X_t = x]=0, \quad \forall  x\in\mathcal X
\label{eq2-2}
\end{align}
Note that $\Psi:\mathcal Q \to  \{ \Psi(\cdot) : \mathcal X \to  [\tau-1, \tau] \}$ is a functional of $q$, $\Psi(q):\mathcal X \to [\tau-1, \tau]$ is a function of $x$  for any fixed $q\in \mathcal Q$, and $\mathcal Y \times  \mathbb R \ni (y,y') \mapsto  \psi_{y}(y') = \tau - \bm 1_{\{y'\le y\}} \in \{\tau-1,\tau\}$. 
Throughout the paper, we use $\| \cdot \|_{\mathcal X}$ as the uniform norm over $\mathcal X$ (i.e., $\|q\|_{\mathcal X}= \sup_{x\in \mathcal X}|q(x)|$), and $\| \cdot \|_{\mathcal Q}$ as the uniform norm over $\mathcal Q$ (i.e., $\|\Psi \|_{\mathcal Q}= \sup_{q\in \mathcal Q} \| \Psi(q) \|_{\mathcal X}=\sup_{q\in \mathcal Q} \sup_{x\in \mathcal X} | \Psi(q)(x) |$).

\begin{remark}
From (\ref{eq2-2}), $\mathbb E[\psi_{q_0(x)}(Y_t)| X_t = x]= \tau - F_{\varepsilon}(q_0(x) -g(x))$ for any $x\in\mathcal X$. 
Since $F_{\varepsilon}$ is a strictly monotonically function from %Assumption
 \ref{A1}, there exists a continuous inverse function denoted by $F_{\varepsilon}^{-1}$. Then, we can write
\begin{align*}
q_0(x)=F_{\varepsilon}^{-1}(\tau)  + g(x),
\end{align*}
and from %Assumption
 \ref{A2}, $q_0$ is uniformly bounded for any $x\in \mathcal X$ (i.e., $q_0\in \mathcal Q$). 
\end{remark}
As the empirical version of the $\Psi$ in (\ref{eq2-2}), we introduce the Generalized Random Forest (GRF) score $\Psi_T$, which will be defined below. 

\subsection{Double-sample tree score and Generalized Random Forest (GRF) score}
Given a vector $\eta=\{Y_0,Y_{-1},\ldots,Y_{1-p}\}$ of initial data independent of $\{\varepsilon_t\}_{t\ge1}$, we suppose that $T$ observations $Y_1,\ldots, Y_T$ from the model (\ref{eq2-1}) are available and we group them in input-output pairs, 
$$\mathcal{D}_T=\left\{(X_1,Y_1),\ldots,(X_T,Y_T)\right\}.$$
In this paper, we construct an conditonal quantile estimator based on the method of Athey et al. (2019). 
 %\cite{athey2019generalized}
 To do so that, we first intoduce a family of subsambles of $\mathcal{D}_T$ by $\mathcal{D}_T\times A \mapsto (\mathcal{I}_s,\mathcal{J}_s)$ where $A$ is an index subset of $\{1,\ldots, T\}$ defined below. 
\begin{definition}\label{def1}
(Double-sample) Let $s=s(T)$ be subsample size and let an index set be $A\subset \left\{1,2,\ldots, T \right\}$ with $|A|=s$. 
A family of the index set $A$ denoted by $\mathcal A_{s}$ is defined as follows:
\begin{align*}
\mathcal A_{s} := \biggl\{ &A = \{A^{\mathcal I},A^{\mathcal J}\},\  A^{\mathcal I},A^{\mathcal J}\subset \left\{1,2,\ldots, T \right\} \biggr| A^{\mathcal{I}}\cap A^{\mathcal{J}}=\emptyset,
    \left|A^{\mathcal{I}}\right| = \left\lfloor \frac{s}{2} \right\rfloor,\  
  \left|A^{\mathcal{J}}\right| = \left\lceil \frac{s}{2} \right\rceil  \biggr\}
\end{align*}
where the elements of $\mathcal A_{s}$ are different from each other. 
In addition, for any $A = \{A^{\mathcal I},A^{\mathcal J}\} \in \mathcal A_{s}$, subsamples $\mathcal{I}_{s}$ and $\mathcal{J}_{s}$ of $\mathcal{D}_T$ are  defined by $\mathcal{I}_{s} = \mathcal{D}_{A^{\mathcal I}}$ and $\mathcal{J}_{s} = \mathcal{D}_{A^{\mathcal J}}$ with $\mathcal{D}_{A^{\cdot}} = \{(X_t,Y_t)\}_{t\in A^{\cdot}}$, respectively. 
\end{definition}
In the double-smple tree defined below, achieves ``honesty'' by dividing its training subsamples into two halves $\mathcal{J}_s$ and $\mathcal{J}_s$. Then, $\mathcal{J}_s$-sample is used to place the splits, while holding out the $\mathcal{I}_s$-sample to do within-leaf estimation falling in each leaf. 
In what follows, we define the splitting rule only by using $\mathcal{J}_s$-sample. 
\begin{definition}\label{def2}
(Splitting rule) Given $\mathcal{J}_s$-sample, we define a sequence of partitions $\mathcal P_0,\mathcal{P}_1,\cdots$ by starting from $\mathcal P_1=\{\mathbb R^p\}$ and then, for each $n\ge1$, construct $\mathcal P_{n+1}$ from $\mathcal P_{n}$ by replacing one set (parent node) $P\in\mathcal P_{n}$ by (childe node) $C_1:=\{x=(x_1,\ldots, x_p) \in P | x_{\xi} \le\zeta\}$ and $C_2:=\{x=(x_1,\ldots, x_p) \in P |x_{\xi}>\zeta\}$, where the split direction $\xi \in \{1,\ldots, p\}$ is randomly chosen 
\footnote{In practice, the optimal direction is chosen from $mtry$ directions at each step of the division where $mtry\sim\min\{\max\{\mathrm{Poisson}(m),1\},p\}$ for some $m\in \mathbb N$.}, and the split position $\zeta =\zeta(\xi) \in \{Y_{t-\xi} | X_t=(Y_{t-1},\ldots, Y_{t-p}) \in P \}$ is chosen to maximize a criterion $\Delta (C_1,C_2)$.
\end{definition}
In this paper, the criterion $\Delta(C_1,C_2)$ is the same as that of Athey et al. (2019)
%\cite{athey2019generalized}
. 
Furthermore, we impose the following assumptions for the splitting rule. 
\begin{assumption}
\begin{enumerate}
\setcounter{enumi}{2}
    \item ($\omega$-Regular)\label{A3} Every split puts at least  a fraction $\omega$ of the observations (of $\mathcal J_s$-sample) in the parent node into each child node, with $\omega\in(0,0.2]$. 
    \item (Random Split)\label{A4} At every split, the probability that the tree splits on the $j$-th feature (i.e., $(X_t)_j = Y_{t-j}$) is bounded from below by some $\pi>0$, for all $j=1,\ldots, p$. 
    \item (PNN (Potential Nearest Neighbor) k-set)\label{A5} There are between $k$ and $2k-1$ observations (of $\mathcal I_s$-sample) in each terminal node . 
    \item (Subsample Size) \label{A6} Subsample size $s$ scales $s=T^{\beta}$ for some $\beta_{\mathrm{min}}<\beta<1$ with 
    $$ \beta_{\mathrm{min}} = 1-\left( 1+ \frac{1}{\pi} \left( \log\left(\omega^{-1}\right)\right)  / \left( \log\left( \left(1-\omega\right)^{-1}\right)\right)\right)^{-1}.$$ 
\end{enumerate}
\end{assumption}
\begin{remark}
This assumption on $\beta$ is the same as (13) of Athey et al. (2019)
%\cite{athey2019generalized}
, in which $s/T\to0$ and $ s\to\infty$ as $T\to\infty$ is satisfied. 
\end{remark}
Under this splitting rule, we denote a given partition of the feature space $\mathbb R^p$ by $\Lambda$, and the subspace (leaf) of rectangular type created by the partitioning by $L_{\ell}$ ($\ell =1,\ldots,|\Lambda|$). 
Then, 
\[ \Lambda=\Lambda(\mathcal I_s^X, \mathcal{J}_s; \xi)=\{L_1,\ldots, L_{|\Lambda|}\},\ \ \mathbb R^p =\bigotimes_{\ell=1}^{|\Lambda|} L_{\ell},\ \  L_{\ell}\cap  L_{\ell'}=\varnothing\ (\ell\neq \ell') \]
where $\mathcal{I}_s^X:= \{X_t=(Y_{t-1},\ldots, Y_{t-p})\}_{t\in A^{\mathcal{I}}}$, $\xi=(\xi_1,\ldots,\xi_{|\Lambda|})$ with $\xi_{\ell}$  ($\ell =1,\ldots,|\Lambda|$) being the split direction for $L_{\ell}$ satisfying %Assumption 
\ref{A4}, independently, each other. 

In addition, we introduce a map $\iota_h$ which transforms from the input vector $ X_t=(Y_{t-1},\ldots,Y_{t-p}) \in \mathbb R^p$ to $[0,1]^p$ by 
\begin{align*}
\iota_h : (y_1,\ldots, y_p) \mapsto (F_h(y_1),\ldots, F_h(y_p))
\end{align*}
where $F_h$ is a cumulative distribution function defined by $F_h(y)=\int_{\infty}^{y}h(\tilde y)d\tilde y$ with 
\begin{align}
h(y):= \frac{1-\bar{\zeta}^{-1}}{\bar{\zeta}-\bar{\zeta}^{-1}}f_{\epsilon}(y+M)+\frac{\bar{\zeta}-1}{\bar{\zeta}-\bar{\zeta}^{-1}}f_{\epsilon}(y-M) ,\label{eq2-3}
\end{align}
$\bar{\zeta}:= \sup_{y\in \mathbb R}\frac{F_{\epsilon}(y+M)}{F_{\epsilon}(y-M)} \in (1,\infty)$ and $M=\sup_{x\in \mathbb R^p}g(x)$.
Based on the map $\iota_h$, we define $Z_t:= \iota_h( X_t)$.
Let $\tilde \Lambda$ be a partition of feature space $[0,1]^p$ in the same manner of $\Lambda$ with $\mathcal I_s^X$ and $\mathcal J_s^X=\{X_t\}_{t\in A^{\mathcal{J}}}$ replaced by $\mathcal {I}_s^Z: =\{Z_t\}_{t\in A^{\mathcal I}}$ and $\mathcal {J}_s^Z: =\{Z_t\}_{t\in A^{\mathcal J}}$. 
Also, let $\tilde L(z)\in \tilde \Lambda$ be the leaf containing the test point $z:= \iota_h(x)$ transformed from the (original) test point $x\in \mathcal X$ into $[0,1]^p$.  
Then, the moment bound of $\mathrm{diam}(\tilde L(z))$ holds (see Lemma \ref{lem1}). 
Moreover, the same moment bound of $\mathrm{diam}(L(x))$ also holds (see Corollary \ref{cor1}). 

By $\Lambda=\Lambda(\mathcal I_s^X, \mathcal J_s; \xi)$ and $\mathcal{I}_s^Y :=\{Y_t\}_{t\in A^{\mathcal{I}}}$, the double-sample tree score is defined as follows. 
\begin{definition}\label{def3}
(Double-sample tree score) Under an observed data $\mathcal{D}_T$, a random vector $\xi$, and any fixed $A\in \mathcal{A}_s$, we fix a partition $\Lambda$ of the feature space $\mathbb R^p$ by Definition \ref{def2}. 
Then, for any $q\in \mathcal Q$ and $x\in \mathcal X$, the double-sample tree score $\mathcal{T}$ is defined by
\[ \mathcal{T}(q ; \mathcal{I}_s, \mathcal J_s, \xi)(x) = \sum_{t\in A^{\mathcal{I}}} \frac{\bm1_{\{ X_t \in L(x)\}}}{\sharp L(x)}  \psi_{q(x)} (Y_t) \]
where $L(x)\in \Lambda(\mathcal I_s^X, \mathcal J_s;\xi)$ is a leaf containing the test point $x$; $\sharp L(x) = |\{t\in A^{\mathcal{I}}: X_t\in L(x)\}|$; and $\psi_{q(x)}(Y_t)=\tau-\bm1_{\{Y_t\le q(x)\}}$.
\end{definition}
%\subsection{Generalized Random Forest (GRF) score}
Gathering the double-sample tree scores, we introduce the GRF score. 
\begin{definition}\label{def4}
(Generalized Random Forest (GRF) score) Let $\mathcal{T}$ be the double-sample tree score by Definition \ref{def3}. Then, for any $q\in \mathcal Q$ and $x\in \mathcal X$, the GRF score $\Psi_T$ is defined by
\begin{align*}
\Psi_T(q)(x) &=  \frac{1}{|\mathcal A_{s}|} \sum_{A\in \mathcal A_{s}}  \mathcal{T}(q ; \mathcal{I}_s, \mathcal J_s, \xi)(x) 
 = \sum_{t=1}^T  \alpha_t(x) \psi_{q(x)} (Y_t)
\end{align*}
where 
\begin{align*}
\alpha_t(x) = \frac{1}{|\mathcal A_{s}|} \sum_{A\in \mathcal A_{s}} \alpha_{A,t}(x), \quad \alpha_{A,t}(x)= \bm 1_{\{t\in A^{\mathcal I}\}}\frac{\bm1_{\{ X_t \in L(x)\}}}{\sharp L(x)} 
\end{align*}
and
\begin{align*}
|\mathcal A_{s}|= \left(
  \begin{array}{cc}
    \multicolumn{2}{c}{T}  \\
    \left\lfloor s/2 \right\rfloor   & \left\lceil s/2 \right\rceil   \\
  \end{array}
\right)
=\left(
  \begin{array}{c}
    T  \\
    s   \\
  \end{array}
\right)\left(
  \begin{array}{c}
    s \\
    \left\lfloor s/2 \right\rfloor \\
  \end{array}
\right)
=\frac{T!}{\left\lfloor s/2 \right\rfloor!\left\lceil s/2 \right\rceil!(T-s)!}.
\end{align*}
\end{definition}

\subsection{Time Series Quantile Regression estimator}
Based on $\Psi_T$ defined by Definition \ref{def4}, we intoroduce a conditional quantile function estimator $\hat{q}_T$ as follows. 
\begin{definition}\label{def5}
Under an observed data $\mathcal{D}_T$, the GRF score $\Psi_T$ for any $x \in \mathcal X$ by Definition \ref{def4}. 
Then, the conditional quantil function estimator $\hat{q}_T=(\hat{q}_T(x))_{x \in \mathcal X}$ is defined by 
\begin{align}
\hat{q}_T( x) 
= \inf \left\{ y\in \mathbb R : \sum_{t=1}^T \alpha_t(x)\left( \tau - \bm 1_{\{Y_t\le y\}} \right) \le 0 \right\}.
\end{align}
\end{definition}
Here we impose the following assumption in order to guarantee the consistency of our estimator. 
\begin{assumption}
\begin{enumerate}
\setcounter{enumi}{6}
\item \label{A7} There exsists some constant $C> 0$, such that, for almost surely, 
$$\left\| \sum_{t=1}^T\alpha_t(\cdot) \psi_{\hat{q}_T(\cdot)}(Y_t) \right\|_{\mathcal X} \le C \max_{t\in \{1,\ldots, T\}} \left\| \alpha_t(\cdot) \right\|_{\mathcal X}.$$
\end{enumerate}
\end{assumption}
This assumption corersponds to the assumption 5 of Athey et al. (2019)
%\cite{athey2019generalized}
. 
Then, we have our main result, that is, uniformly consistency of the conditionl quantile function estimator $\hat{q}_T$. 
\begin{theorem}\label{thm1}
Under %Assumptions
 \ref{A1} -  \ref{A7}, the conditional quantile function estimator $\hat{q}_T$ converges to the true conditional quantile function $q_0$, in probability, that is, 
\begin{align*}
\|\hat{q}_T - q_0 \|_{\mathcal X}\stackrel{p}{\to} 0 \quad \mathrm{as}\ T\to\infty.
\end{align*}
\end{theorem}
The proofs of the theorems, lemmas and corollary are given in Appendix \ref{app}.\\

When $T$ is large, it is not realistic to generate all possible $|\mathcal A_s|$ types of tree. 
In practice, under a sufficiently large $B$, we randomly choose $B$ types of subset of $\{1,\ldots, T\}$ defined by $\{A^{(b)}\in \mathcal A_s\}_{b=1,\ldots, B}$, and generate a double-sample tree score $\mathcal T^{(b)}$ based on the subsample $(\mathcal I_s^{(b)},\mathcal J_s^{(b)})$ for each $A^{(b)}$.  
\begin{definition}\label{def6}
%\paragraph{\textbf{Definition 4'}}
%\textit{
(Generalized Random Forest (GRF) score) Define $\mathcal A_s^{B}=\{A^{(b)} = \{A_b^{\mathcal I}, A_b^{\mathcal J}\} \in \mathcal A_s\}_{b=1,\ldots, B}$. 
For each subsamples $\{\mathcal I_s^{(b)}, \mathcal J_s^{(b)}\}$ determined by $\mathcal D_T$ and $A^{(b)}\in \mathcal A_s^B$, we generate the double-sample tree score $\mathcal{T}$ by Definition \ref{def3}. 
Then, for any $q\in \mathcal Q$ and $x\in \mathcal X$, the GRF score $\Psi_T^{B}$ is defined by
\begin{align*}
\Psi_T^{B}(q)(x) =  \frac{1}{B} \sum_{b=1}^B \mathcal{T}(q ; \mathcal{I}_s^{(b)}, \mathcal J_s^{(b)}, \xi)(x) 
 = \sum_{t=1}^T  \alpha_t^{B}(x) \psi_{q(x)} (Y_t)
\end{align*}
where 
\begin{align*} 
\alpha_t^{B}(x) = \frac{1}{B} \sum_{b=1}^B \alpha_{A^{(b)},t}(x),\quad  \alpha_{A^{(b)},t}(x) = \bm 1_{\{t\in A_b^{\mathcal I }\}}\frac{\bm1_{\{ X_t \in L^{(b)}(x)\}}}{\sharp L^{(b)}(x)} 
\end{align*}
and $L^{(b)}(x)\in \Lambda (\mathcal I_s^{X (b)} , \mathcal J_s^{(b)};\xi)$ is a leaf containing the test point $x$; $\sharp L^{(b)}(x)= |\{t\in A_b^{\mathcal I}:X_t \in L^{(b)}(x)\}|$.
%}
\end{definition}

Based on $\Psi_T^B$ defined by Definition \ref{def6}%4'
, we introduce another conditional quantile function estimator $\hat{q}_T^{(B)}=(\hat{q}_T^{(B)}(x))_{x \in \mathcal X}$ as followns. 
\begin{definition}\label{def7}
%\paragraph{\textbf{Definition 5'}}
%\textit{
Under an observed data $\mathcal{D}_T$, the GRF score $\Psi_T^B$ for all $x\in \mathcal X$ by Definition \ref{def6}%4'
. 
Then, the conditional quantile function estimator $\hat{q}_T^B$ is defined by
\begin{align}
\hat{q}_T^B( x) 
= \inf \left\{ y\in \mathbb R : \sum_{t=1}^T \alpha_t^B(x)\left( \tau - \bm 1_{\{Y_t\le y\}} \right) \le 0 \right\}.
\end{align}
%}
\end{definition}

If $B$ is sufficiently large, we have uniformaly consistency of $\hat{q}_T^B$ as follows. 
\begin{theorem}\label{thm2}
Under %Assumptions 
\ref{A1} - \ref{A7}, and $B\equiv B(T)$ with $\lim_{T\to \infty} \frac{1}{B} =0$, $\hat{q}_T^B$ converges to $q_0$, in probability, that is,
\begin{align*}
\|\hat{q}_T^B - q_0 \|_{\mathcal X} 
\stackrel{p}{\to} 0 \quad \mathrm{as}\ T\to\infty.
\end{align*}
\end{theorem}

\section{Simulation}\label{sec3}
In this section, we first check the asymptotic properties of tsQRF for several data generation settings. 
 To illustrate the characteristics of our method, we compared tsQRF and WNW estimators (Cai, 2002)
 %\cite{cai2002regression})
 . 

\subsection{Data generation process}\label{sec3.1}

Four data generating models were used in this simulation. 
First, model (a) generates data using a bounded oscillating function, which was used by Davis and Nielsen (2020)%\cite{davis2020modeling}
.

\renewcommand{\theenumi}{(\alph{enumi})}
\begin{enumerate}
%    \item First-order Markov chain model（Davis and Nielsen, %2020;
%     equation (4.1) %) 
%     \cite{davis2020modeling}) \label{f}
    \item First order Markov chain model (Davis and Nielsen, 2020
    %\cite{davis2020modeling}
    , equation (4.1) ) \label{f}
    \begin{align*}
        Y_t=\cos(5Y_{t-1})e^{-Y_{t-1}^2}+\varepsilon_t
    \end{align*}
\end{enumerate}

This model satisfies the boundedness assumption \ref{A2} for function $g$. 
However, such a bounded model is often not used in an actual time series data analysis. 
Therefore, we also generate data from models that do not satisfy %assumption 
\ref{A2}, such as models \ref{AR2}, \ref{TAR2}, and \ref{AR5}, and examine how our estimator converges to the true quantile function.

\begin{enumerate}
\setcounter{enumi}{1}
    \item AR(2) model \label{AR2}
    \begin{align*}
    Y_t = 0.5Y_{t-1} + 0.4Y_{t-2} + \varepsilon_t
    \end{align*}
    \item Non-linear AR(2) model \label{TAR2}
    \begin{align*}
    Y_t = \left\{\begin{array}{ll}
    2.9 - 0.4Y_{t-1} - 0.1Y_{t-2} + \varepsilon_t  & (\mathrm{if}\ Y_{t-1} \leq 1)\\
     -1.5 + 0.2Y_{t-1} + 0.3Y_{t-2} + \varepsilon_t  &(\mathrm{if}\ Y_{t-1} > 1)
    \end{array}
    \right.
    \end{align*}
    \item AR(5) model \label{AR5}
    \begin{align*}
    Y_t = 0.7Y_{t-1} - 0.6Y_{t-2} + 0.4Y_{t-3} -0.2 Y_{t-4} +0.1 Y_{t-5} + \varepsilon_t
    \end{align*}
\end{enumerate}

For each model, we consider two types of error distribution: the standard normal distribution (normal) and the standard Laplace distribution (Laplace) for $\varepsilon_1,\varepsilon_2,\ldots$. 
Both distributions satisfy %assumption 
\ref{A1}%(A-1)
. 
The difference between the two error distributions is the behavior of the tail of the distribution. 
In the simulation, we examined the effects of the tail behavior of error distributions. 
For each model, the true value of the conditional $\tau$-quantile given $X_t=x_t=(y_{t-1},\ldots,y_{t-p})$ is
\begin{align*}
 q_0(x_t) = g(x_t) + F_{\varepsilon}^{-1}(\tau),
\end{align*}
where $F_{\varepsilon}^{-1}$ is the inverse of the distribution function of the error term.

In the simulation, $T+T'$ length time series data were generated from each model. 
The first $T(:=1000,5000)$ is used as training data for estimating the models, and the remaining $T'(:=500)$ is used as test data to evaluate the accuracy of the quantile prediction. 
We generated $R(:=100)$ replicates of these time series data to compute the estimation/prediction error.

\subsection{Evaluation of estimation accuracy for training data (consistency)}\label{sec3.2}

First, we illustrate the consistency of the conditional quantile estimated using the tsQRF. 
We estimated the $\tau=1\%,10\%,50\%,90\%$, and $99\%$ quantiles for each scenario. 
We used the R package \textbf{grf} to estimate the target quantiles and set the parameters of the GRF as subsample size $s=\frac{T}{2}$, number of trees $B=2000$, and the parameter $\omega=0.05$ for the ratio of splits (these are all default values from the \textbf{grf} package)

Theoretically, $|\mathcal{A}_s|$ is used for $B$; however, in practice, when $T$ and $s$ are sufficiently large, $|\mathcal{A}_s|$ becomes very large, and the computational cost becomes expensive. 
Therefore, instead of $|\mathcal{A}_s|$, we restrict the number of trees to $B\ll|\mathcal{A}_s|$ ($B$ elements of $A$ from $\mathcal{A}_s$ are randomly chosen). 
If $B$ is as large as the sample size, the approximation error can be sufficiently small (Wager and Athey (2018)%\cite{wager2018estimation}
).

In this simulation, we computed the mean and standard deviation of the bias of the estimates and the mean squared error. 
Let $y_{1}^{(r)},\ldots,y_{T}^{(r)}$ be the training dataset for each replicate $r\ (r=1,\ldots,R)$. 
We defined the difference between the true value of the quantiles and the estimated conditional quantiles at $x_{1}^{(r)},\ldots, x_{T}^{(r)}$ as $\mathrm{Bias}_t^{(r)}=\hat q_T(x_{t}^{(r)})-q_0(x_{t}^{(r)})$. 
The average is defined as $\mathrm{Bias}^{(r)}=\frac{1}{T}\sum_{t=1}^{T}\mathrm{Bias}_t^{(r)}$. 
 The mean (MBias), standard deviation (SDBias), and mean squared error (MSE) of $\mathrm{Bias}^{(1)},\ldots,\mathrm{Bias}^{(R)}$ are defined as follows.

\begin{align*}
    \mathrm{MBias}&=\frac{1}{R}\sum_{r=1}^{R}\mathrm{Bias}^{(r)} \\
    \mathrm{SDBias}&=\sqrt{\frac{1}{R-1}\sum_{r=1}^{R}(\mathrm{Bias}^{(r)}-\mathrm{MBias})^2} \\
    \mathrm{MSE}&=\frac{1}{R}\sum_{r=1}^{R}\frac{1}{T}\sum_{t=1}^{T}\mathrm{Bias}_t^{(r)\ 2} 
\end{align*}

We performed the simulation for each data generating model (a) $\sim$ \ref{AR5} with the length of time series $T=1000$ and $T=5000$, and two types of error distribution (normal and Laplace). 
The MBias, SDBias and MSE for each simulation scenario are summarized in Tables \ref{mbias_convergence_train}$\sim$\ref{mse_convergence_train}. 

\begin{table}
\centering
\caption{MBias (training data)}\label{mbias_convergence_train}
\begin{tabular}{cccccccc}
  %\hline
Model &$\varepsilon_t$&&1\% & 10\% & 50\% & 90\% & 99\% \\  
   \hline
   (a)&Normal & $T=1000$ &0.297 & 0.021 & 0.004 & \textcolor{red}{-0.012} & \textcolor{red}{-0.288} \\ 
  &&$T=5000$ &\textcolor{red}{0.297} & \textcolor{red}{0.018} & \textcolor{red}{0.000} & -0.017 & -0.294 \\ 
  &Laplace & $T=1000$ &0.652 & \textcolor{red}{-0.017} & -0.003 & \textcolor{red}{0.015} & -0.667 \\ 
   && $T=5000$ &\textcolor{red}{0.633} & -0.026 & \textcolor{red}{0.000} & 0.025 & \textcolor{red}{-0.645} \\ 
  \hline
  \ref{AR2}&Normal & $T=1000$ & \textcolor{red}{-0.007} & -0.072 & -0.011 & 0.048 & \textcolor{red}{-0.005} \\ 
  &&$T=5000$ & 0.078 & \textcolor{red}{-0.014} & \textcolor{red}{-0.000} & \textcolor{red}{0.014} & -0.074 \\ 
  &Laplace & $T=1000$ &\textcolor{red}{0.168} & -0.133 & -0.010 & 0.105 & \textcolor{red}{-0.185} \\ 
   && $T=5000$ &0.209 & \textcolor{red}{-0.064} & \textcolor{red}{-0.005} & \textcolor{red}{0.051} & -0.217 \\ 
   \hline
  \ref{TAR2}&Normal & $T=1000$ & \textcolor{red}{0.012} & -0.039 & -0.003 & 0.046 & \textcolor{red}{0.057} \\ 
  &&$T=5000$ & 0.077 & \textcolor{red}{-0.007} & \textcolor{red}{0.001} & \textcolor{red}{0.009} & -0.064 \\ 
   &Laplace & $T=1000$ &\textcolor{red}{0.182} & -0.092 & -0.010 & 0.070 & \textcolor{red}{-0.175} \\ 
   && $T=5000$ &0.254 & \textcolor{red}{-0.026} & \textcolor{red}{-0.002} & \textcolor{red}{0.033} & -0.211 \\ 
   \hline
   \ref{AR5}&Normal & $T=1000$ &-0.281 & -0.163 & \textcolor{red}{0.000} & 0.161 & 0.263 \\ 
  &&$T=5000$ &-\textcolor{red}{0.170} & \textcolor{red}{-0.108} & -0.002 & \textcolor{red}{0.107} & \textcolor{red}{0.168} \\ 
  &Laplace & $T=1000$ &-0.260 & -0.289 & \textcolor{red}{0.001} & 0.277 & 0.210 \\ 
   && $T=5000$ &\textcolor{red}{-0.099} & \textcolor{red}{-0.183} & -0.004 & \textcolor{red}{0.176} & \textcolor{red}{0.102} \\ 
   \hline
\end{tabular}
\end{table}

\begin{table}
\centering
\caption{SDBias (traning data)}\label{sdbias_convergence_train}
\begin{tabular}{cccccccc}
  %\hline
Model &$\varepsilon_t$&&1\% & 10\% & 50\% & 90\% & 99\% \\  
  \hline
  (a)&Normal & $T=1000$ &0.068 & 0.046 & 0.037 & 0.054 & 0.081 \\ 
  &&$T=5000$ &\textcolor{red}{0.036} & \textcolor{red}{0.024} & \textcolor{red}{0.015} & \textcolor{red}{0.025} & \textcolor{red}{0.038} \\ 
  &Laplace & $T=1000$ &0.218 & 0.093 & 0.037 & 0.086 & 0.201 \\ 
   && $T=5000$ &\textcolor{red}{0.088} & \textcolor{red}{0.042} & \textcolor{red}{0.016} & \textcolor{red}{0.042} & \textcolor{red}{0.096} \\ 
   \hline
  \ref{AR2}&Normal & $T=1000$ & 0.084 & 0.045 & 0.031 & 0.048 & 0.083 \\ 
  &&$T=5000$ & \textcolor{red}{0.041} & \textcolor{red}{0.019} & \textcolor{red}{0.014} & \textcolor{red}{0.022} & \textcolor{red}{0.039} \\ 
  &Laplace & $T=1000$ &0.222 & 0.089 & 0.039 & 0.081 & 0.260 \\ 
   && $T=5000$ &\textcolor{red}{0.112} & \textcolor{red}{0.038} & \textcolor{red}{0.016} & \textcolor{red}{0.042} & \textcolor{red}{0.102} \\ 
   \hline
  \ref{TAR2}&Normal & $T=1000$ & 0.093 & 0.047 & 0.038 & 0.054 & 0.085 \\ 
  &&$T=5000$ & \textcolor{red}{0.041} & \textcolor{red}{0.020} & \textcolor{red}{0.015} & \textcolor{red}{0.021} & \textcolor{red}{0.037} \\ 
   &Laplace & $T=1000$ &0.254 & 0.093 & 0.035 & 0.086 & 0.240 \\ 
   && $T=5000$ & \textcolor{red}{0.101} & \textcolor{red}{0.040} & \textcolor{red}{0.017} & \textcolor{red}{0.040} & \textcolor{red}{0.111} \\ 
   \hline
   \ref{AR5}&Normal & $T=1000$ &0.109 & 0.054 & 0.034 & 0.048 & 0.105 \\ 
  &&$T=5000$ &\textcolor{red}{0.044} & \textcolor{red}{0.023} & \textcolor{red}{0.015} & \textcolor{red}{0.025} & \textcolor{red}{0.049} \\ 
  &Laplace & $T=1000$ &0.280 & 0.099 & 0.035 & 0.088 & 0.244 \\ 
   && $T=5000$ &\textcolor{red}{0.117} & \textcolor{red}{0.036} & \textcolor{red}{0.016} & \textcolor{red}{0.040} & \textcolor{red}{0.112} \\ 
   \hline
\end{tabular}
\end{table}

\begin{table}
\centering
\caption{MSE (traning data)}\label{mse_convergence_train}
\begin{tabular}{cccccccc}
  %\hline
Mode &$\varepsilon_t$&&1\% & 10\% & 50\% & 90\% & 99\% \\  
  \hline
  (a)&Normal & $T=1000$ &0.316 & 0.118 & \textcolor{red}{0.062} & 0.117 & \textcolor{red}{0.304} \\
  &&$T=5000$ &\textcolor{red}{0.309} & \textcolor{red}{0.116} & 0.063 & \textcolor{red}{0.116} & 0.310 \\ 
  &Laplace & $T=1000$ &1.941 & 0.394 & 0.057 & \textcolor{red}{0.386} & \textcolor{red}{1.872} \\ 
   && $T=5000$ &\textcolor{red}{1.865} & \textcolor{red}{0.393} & \textcolor{red}{0.056} & 0.397 & 1.875 \\ 
   \hline
  \ref{AR2}&Normal & $T=1000$ &0.182 & 0.091 & 0.055 & 0.090 & 0.189 \\ 
  &&$T=5000$ &\textcolor{red}{0.156} & \textcolor{red}{0.067} & \textcolor{red}{0.038} & \textcolor{red}{0.066} & \textcolor{red}{0.158} \\ 
  &Laplace & $T=1000$ &1.211 & 0.287 & 0.085 & 0.273 & 1.195 \\ 
   && $T=5000$ &\textcolor{red}{1.137} & \textcolor{red}{0.225} & \textcolor{red}{0.047} & \textcolor{red}{0.223} & \textcolor{red}{1.175} \\ 
   \hline
  \ref{TAR2}&Normal & $T=1000$ &0.368 & 0.150 & 0.054 & 0.130 & 0.492 \\ 
  &&$T=5000$ & \textcolor{red}{0.207} & \textcolor{red}{0.075} & \textcolor{red}{0.037} & \textcolor{red}{0.071} & \textcolor{red}{0.225} \\ 
   &Laplace & $T=1000$ &1.360 & 0.313 & 0.064 & 0.285 & 1.318 \\ 
   && $T=5000$ & \textcolor{red}{1.163} & \textcolor{red}{0.221} & \textcolor{red}{0.036} & \textcolor{red}{0.219} & \textcolor{red}{1.181} \\ 
   \hline
   \ref{AR5}&Normal & $T=1000$ &0.335 & 0.189 & 0.145 & 0.188 & 0.321 \\ 
  &&$T=5000$ &\textcolor{red}{0.183} & \textcolor{red}{0.108} & \textcolor{red}{0.084} & \textcolor{red}{0.108} & \textcolor{red}{0.179} \\ 
  &Laplace & $T=1000$ &1.155 & 0.471 & 0.299 & 0.469 & 1.058 \\ 
   && $T=5000$ &\textcolor{red}{0.794} & \textcolor{red}{0.281} & \textcolor{red}{0.161} & \textcolor{red}{0.279} & \textcolor{red}{0.791} \\ 
   \hline
\end{tabular}
\end{table}

\color{black}

We first discuss the case where $\varepsilon_t$ follows a standard normal distribution. 
In model (a), there is no significant difference between $T=1000$ and $T=5000$ for MBias and MSE, but SDBias is much closer to zero when $T=5000$. 
In model \ref{AR2}, SDBias and MSE are much closer to zero when $T=5000$ than $T=1000$. 
For MBias, $T=5000$ is closer to zero for $\tau=10\%,50\%,90\%$, but for relatively high (low) levels, such as $\tau=1\%,99\%$, $T=1000$ is closer to zero. 
Similarly in model \ref{TAR2}, SDBias and MSE are much closer to zero when $T=5000$, and for MBias, $T=5000$ is closer to zero for $\tau=10\%,50\%,90\%$ compared to $T=1000$.
Because models \ref{AR2} and \ref{TAR2} have the same lag order (i.e., dimension of covariate space), there is no significant difference in any of the indices, except for the relatively high (low) level. 
When the lag order is high, as in model \ref{AR5}, MBias is not significantly different, although $T=1000$ is closer to zero at $\tau=50\%$, and $T=5000$ is closer to zero  at all other levels. 
The accuracy of the estimation decreases slightly as the lag order increases, and in most of the other levels, the model takes values farther from zero than the other models. 

Next, we discuss the case in which $\varepsilon_t$ follows a standard Laplace distribution. 
In model (a), as in the case of the standard normal distribution, there is no significant difference between $T=1000$ and $T=5000$ for MBias and MSE, whereas for SDBias, the value approaches zero when $T=5000$. At relatively high (low) levels, such as $\tau=1\%$ and $\tau=99\%$, the effect of increasing $T$ is smaller for the MBias and MSE than for the other models.
In models \ref{AR2} and \ref{TAR2}, as in the case of the standard normal distribution, SDBias and MSE are much closer to zero when $T=5000$, and for MBias, $T=5000$ is closer to zero for $\tau=10\%, 50\%, 90\%$. 
For MBias, $T=5000$ is closer to zero than $T=10\%, 50\%, 90\%$. 
As in the case of the standard normal distribution, $T=5000$ is closer to zero at most levels in model \ref{AR5}, and is farther from zero than the other models.

Comparing the case of $\varepsilon_t$ with the standard normal distribution and the standard Laplace distribution, MBias, SDBias and MSE are farther from zero for the standard Laplace distribution which has a relatively heavy tail distribution. 
At $\tau=1\%,99\%$, this phenomenon is remarkable. 
Although there are model-specific differences, for most levels, MBias, SDBias, and MSE tend to approach zero as $T$ increases, we may conclude that the quantile estimators by tsQRF are consistent with the simulations. 

\begin{figure}
    \centering
    \includegraphics[height=8.0cm,width=12cm]{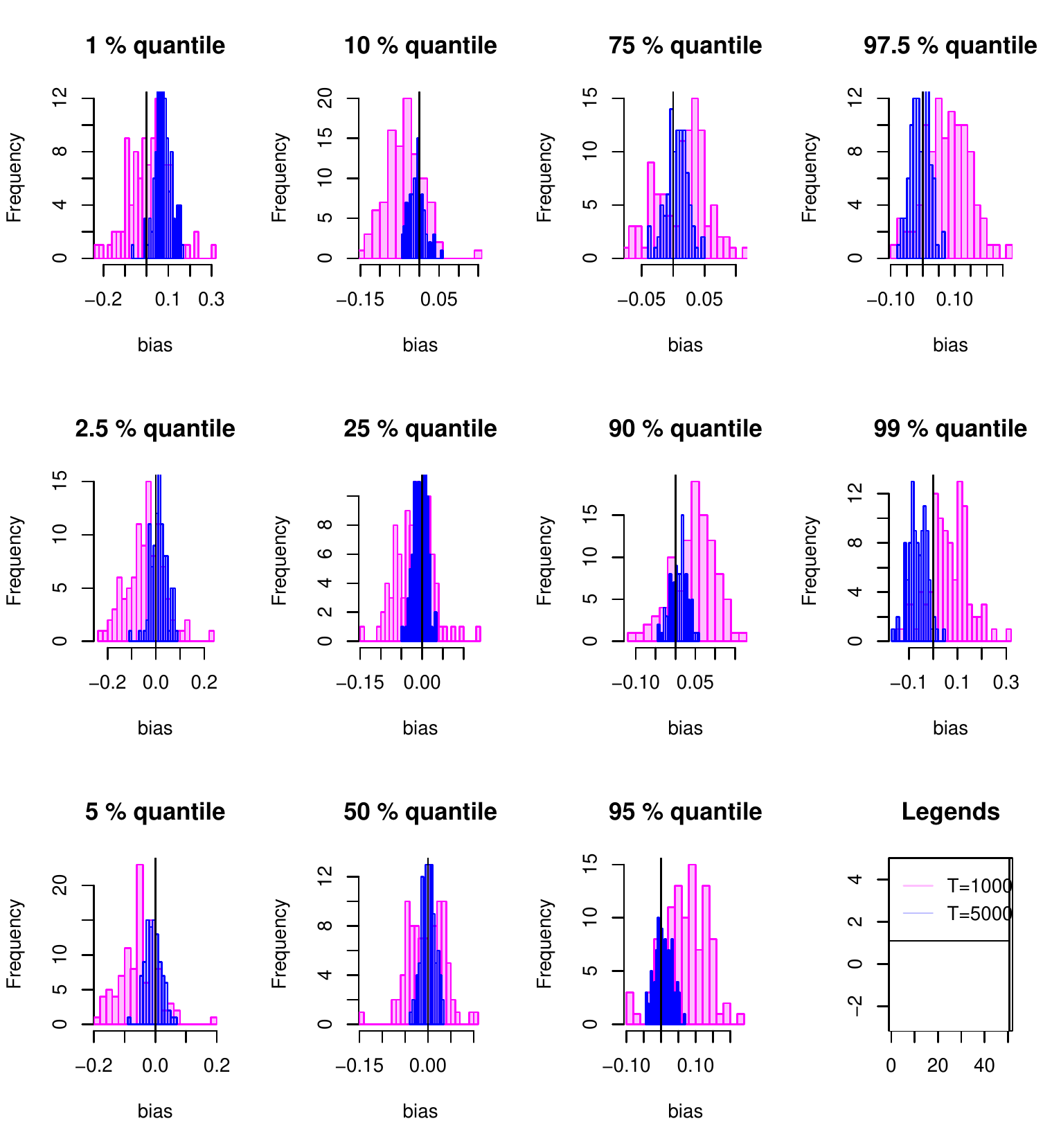}
    \caption{Convergence of estimator for model\ref{TAR2}}
    \label{conv_TAR2_insample}
\end{figure}

Figure \ref{conv_TAR2_insample} shows the histgram of $\mathrm{Bias}^{(1)},\ldots,\mathrm{Bias}^{(R)}$ for model \ref{TAR2} at each $\tau$. 

The consistency of the estimators is illustrated in Figure \ref{conv_TAR2_insample}. 
The variance becomes small, and the bias approaches zero for $T=5000$ compared with $T=1000$ except for $\tau=1\%$ and $99\%$. 
In addition, for $T=5000$, the bias tends to be positive, that is, the estimated value is larger than the true value, as $T=5000$ and $\tau$ approaches $1\%$, the estimated value tends to be smaller than the true value.
From Table \ref{mbias_convergence_train}, this trend can also be observed in models (a) and \ref{AR2}.
The opposite trend was observed for model \ref{AR5}. 
For the previously mentioned features at the 1\% and 99\% quantiles, the shape of the distribution of $\varepsilon_t$ or the boundedness assumption of the function $g$ is considered to be affected.

\subsection{Comparison of tsQRF and other quantile regression methods}\label{sec3.3}

Here, we compare the prediction accuracy when applying the kernel quantile regression of Cai (2002)
%\cite{cai2002regression}
 and  tsQRF proposed in this study.

Cai (2002) %\cite{cai2002regression}
 proposed a nonparametric method for estimating conditional quantiles by determing the inverse function of the Weighted Nadaraya-Watson (WNW) estimator of the conditional distribution function. 
 For strongly stationary and $\alpha$-mixing $\{(Y_t,X_t)\}_{t\geq 1}$, the WNW estimator of the conditional distribution function is defined as follows.

\begin{align*}
    \hat F(y|x)=\frac{\sum_{t=1}^Tp_t(x)K\left(\frac{x-X_t}{h}\right)\bm1_{\{Y_t\le y\}}}{\sum_{t=1}^Tp_t(x)K\left(\frac{x-X_t}{h}\right)}
\end{align*}
Here, $p_t(x)$ is a nonnegative weighting function that, satisfies $\sum_{t=1}^T p_t(x)=1$. 
In this study, we set $p_t(x)=1/T$ using the R package \textbf{np}. 
To reduce the computational cost, we set the parameters of the \textbf{np} package as itmax=5000, tol=0.1, and ftol=0.1. 
In this simulation, we used training data $T=1000$ and the accuracy of the estimators is compared using test data $T' = 500$. 
The quantile levels to be compared were $\tau=10\%,50\%,90\%$, and the prediction accuracy was evaluated using the same measures as in Section 3.2, but $\{(x_t^{(r)}, y_t^{(r)})\}_{t=1}^T$ is replaced by $\{(x_t^{(r)}, y_t^{(r)})\}_{t=T+1}^{T+T'}$ except for $\hat{q}_T$. 
We denote these by $\widetilde{\mathrm{MBias}}$, $\widetilde{\mathrm{SDBias}}$ and $\widetilde{\mathrm{MSE}}$, respectively. 

Tables \ref{bias_test} and \ref{mse_test} show the $\widetilde{\mathrm{MBias}}$, $\widetilde{\mathrm{SDBias}}$ and $\widetilde{\mathrm{MSE}}$ of tsQRF and kernel QR (WNW) for data generation models \ref{f} - \ref{AR5}, respectively.

\begin{table}
\centering
\caption{Comparison of bias for estimator on test data}\label{bias_test}
\begin{tabular}{ccccccccc}
  \hline
  &&&\multicolumn{3}{c}{$\widetilde{\mathrm{MBias}}$}&\multicolumn{3}{c}{$\widetilde{\mathrm{SDBias}}$} \\
  \cmidrule(lr){4-6} \cmidrule(lr){7-9}
 Models &$\varepsilon_t$&&10\% &50\% & 90\% &10\% &50\% & 90\%  \\  
  \hline
  (a)&Normal&WNW &-0.104 & \textcolor{red}{0.003} & 0.111&0.052 & \textcolor{red}{0.035} & 0.058 \\ 
   %&&CAViaR &  &&&&&  \\
  &&tsQRF &  \textcolor{red}{0.019} & 0.004 & \textcolor{red}{-0.008}&  \textcolor{red}{0.051} & 0.039 & \textcolor{red}{0.058}  \\ 
  &Laplace&WNW&-0.107 & \textcolor{red}{-0.001} & 0.112 &\textcolor{red}{0.098} & 0.042 & \textcolor{red}{0.086}  \\
  %&&CAViaR &  &&&&&  \\
  &&tsQRF &\textcolor{red}{-0.015} & -0.002 & \textcolor{red}{0.018}&0.101 & \textcolor{red}{0.041} & 0.097  \\
   \hline
  \ref{AR2}&Normal&WNW &-0.112 & -0.008 & 0.090&0.063 & 0.055 & 0.071  \\ 
  %&&CAViaR &  &&&&&  \\
  &&tsQRF &  \textcolor{red}{-0.078} & \textcolor{red}{-0.005} & \textcolor{red}{0.067}& \textcolor{red}{0.063} & \textcolor{red}{0.052} & \textcolor{red}{0.066}  \\ 
  &Laplace&WNW&\textcolor{red}{-0.127} & \textcolor{red}{0.008} & \textcolor{red}{0.131}&0.133 & 0.091 & \textcolor{red}{0.113} \\ 
  %&&CAViaR &  &&&&&  \\
  &&tsQRF &-0.129 & 0.014 & 0.153&\textcolor{red}{0.130} & \textcolor{red}{0.089} & 0.114  \\
   \hline
   \ref{TAR2}&Normal&WNW &-0.136 & \textcolor{red}{-0.002} & 0.138&0.054 & 0.046 & 0.059  \\ 
   %&&CAViaR &  &&&&&  \\
  &&tsQRF &\textcolor{red}{-0.045}&-0.008&\textcolor{red}{0.048}&\textcolor{red}{0.049}&\textcolor{red}{0.041}&\textcolor{red}{0.057}  \\ 
  &Laplace&WNW&-0.173 & \textcolor{red}{-0.015} & 0.130&0.099 & 0.044 & \textcolor{red}{0.088}  \\
  %&&CAViaR &  &&&&&  \\
  &&tsQRF &\textcolor{red}{-0.100} & -0.018 &\textcolor{red}{0.069}&\textcolor{red}{0.096} & \textcolor{red}{0.038} & 0.090  \\
   \hline
   \ref{AR5}&Normal&WNW &\textcolor{red}{-0.146} & \textcolor{red}{0.002} & \textcolor{red}{0.153}&\textcolor{red}{0.058} & \textcolor{red}{0.040} & \textcolor{red}{0.054} \\ 
   %&&CAViaR &  &&&&&  \\
  &&tsQRF &  -0.180 & 0.003 & 0.184&  0.062 & 0.042 & 0.057  \\ 
  &Laplace&WNW&\textcolor{red}{-0.187} & \textcolor{red}{0.004} & \textcolor{red}{0.190}&0.106 & \textcolor{red}{0.048} & \textcolor{red}{0.087} \\ 
  %&&CAViaR &  &&&&&  \\
  &&tsQRF &-0.296 & 0.005 & 0.292&\textcolor{red}{0.102} & 0.049 & 0.095  \\
   \hline
\end{tabular}
\end{table}

\begin{table}
\centering
\caption{Comparison of $\widetilde{\mathrm{MSE}}$ (test data)}\label{mse_test}
\begin{tabular}{ccccccccc}
  \hline
  %&&&\multicolumn{3}{c}{学習データ}&\multicolumn{3}{c}{テストデータ} \\
  % \cmidrule(lr){4-6} \cmidrule(lr){7-9}
  Model &$\varepsilon_t$& &10\% &50\% & 90\%  \\  
  \hline
  (a)&Normal&WNW &\textcolor{red}{0.059} & \textcolor{red}{0.030} & \textcolor{red}{0.059} \\ 
   %&&CAViaR &  &&  \\
  &&tsQRF & 0.119 & 0.062 & 0.117  \\
  &Laplace&WNW &\textcolor{red}{0.198} & 0.076 & \textcolor{red}{0.191} \\ 
   %&&CAViaR &  &&  \\
  &&tsQRF & 0.400 & \textcolor{red}{0.057} & 0.397  \\
   \hline
  \ref{AR2}&Normal&WNW&\textcolor{red}{0.078} & \textcolor{red}{0.054} & \textcolor{red}{0.077} \\ 
  %&&CAViaR &  &&  \\
  &&tsQRF & 0.109 & 0.075 & 0.109  \\
  &Laplace&WNW&\textcolor{red}{0.307} & 0.175 & \textcolor{red}{0.262} \\ 
  %&&CAViaR &  &&  \\
  &&tsQRF & 0.348 & \textcolor{red}{0.162} & 0.352  \\
   \hline
   \ref{TAR2}&Normal&WNW &0.181 & 0.085 & 0.161  \\
   %&&CAViaR &  &&  \\
  &&tsQRF &\textcolor{red}{0.161} & \textcolor{red}{0.063} & \textcolor{red}{0.135}  \\
  &Laplace&WNW &0.369 & 0.128 & 0.329  \\
   %&&CAViaR &  &&  \\
  &&tsQRF &\textcolor{red}{0.310} & \textcolor{red}{0.078} & \textcolor{red}{0.273}  \\
   \hline
   \ref{AR5}&Normal&WNW &\textcolor{red}{0.151} & \textcolor{red}{0.110} & \textcolor{red}{0.153} \\ 
   %&&CAViaR &  &&  \\
  &&tsQRF & 0.201 & 0.160 & 0.203  \\
  &Laplace&WNW &\textcolor{red}{0.413} & \textcolor{red}{0.256} & \textcolor{red}{0.401} \\ 
   %&&CAViaR &  &&  \\
  &&tsQRF & 0.490 & 0.350 & 0.491  \\
   \hline
\end{tabular}
\end{table}

We first discuss the case where $\varepsilon_t$ follows a standard normal distribution. 
In models (a) and \ref{AR2}, there is no significant difference for $\widetilde{\mathrm{SDBias}}$, while $\widetilde{\mathrm{MBias}}$ is closer to zero for tsQRF, and $\widetilde{\mathrm{MSE}}$ is closer to zero for WNW. 
In model \ref{TAR2}, tsQRF tends to be closer to zero for all indicators. 
Finally, the WNW was closer to zero for all indicators in model \ref{AR5}.

Next, we discuss the case in which $\varepsilon_t$ follows a standard Laplace distribution. 
As in the case of the standard normal distribution, there is no significant difference in $\widetilde{\mathrm{SDBias}}$ in model (a), and $\widetilde{\mathrm{MBias}}$ tends to be closer to zero for tsQRF, whereas $\widetilde{\mathrm{MBias}}$ and $\mathrm{MSE}$ tend to be closer to zero for WNW. 
In model \ref{AR2}, $\widetilde{\mathrm{SDBias}}$ is not significantly different, but $\widetilde{\mathrm{MBias}}$ and $\widetilde{\mathrm{MSE}}$ are closer to zero for WNW. 
In model \ref{TAR2}, as in the case of the standard normal distribution, tsQRF was approximately zero for all indicators. 
Finally, model \ref{AR5} shows that the WNW is approximately zero for all indicators.

\section{Empirical Results}\label{sec4}
\begin{figure}[htbp]
        \begin{minipage}[b]{0.46\textwidth}
            \centering
            \includegraphics[width=\columnwidth]{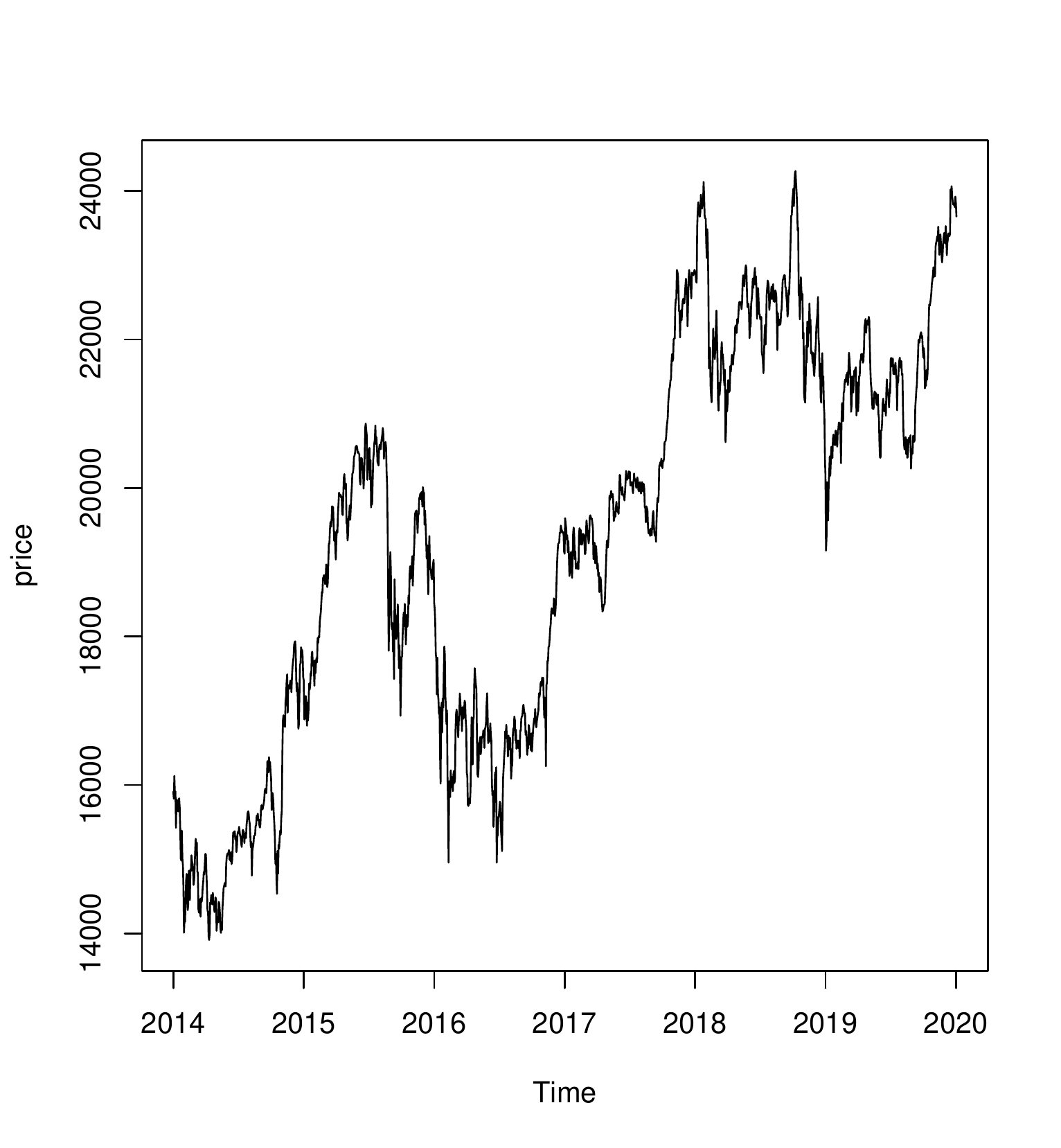}
            \caption{Time series of price $p_t$}\label{plot_nikkei}
        \end{minipage}
        \hspace{0.04\textwidth} % ここで隙間作成
        \begin{minipage}[b]{0.46\textwidth}
            \centering
            \includegraphics[ width=\columnwidth]{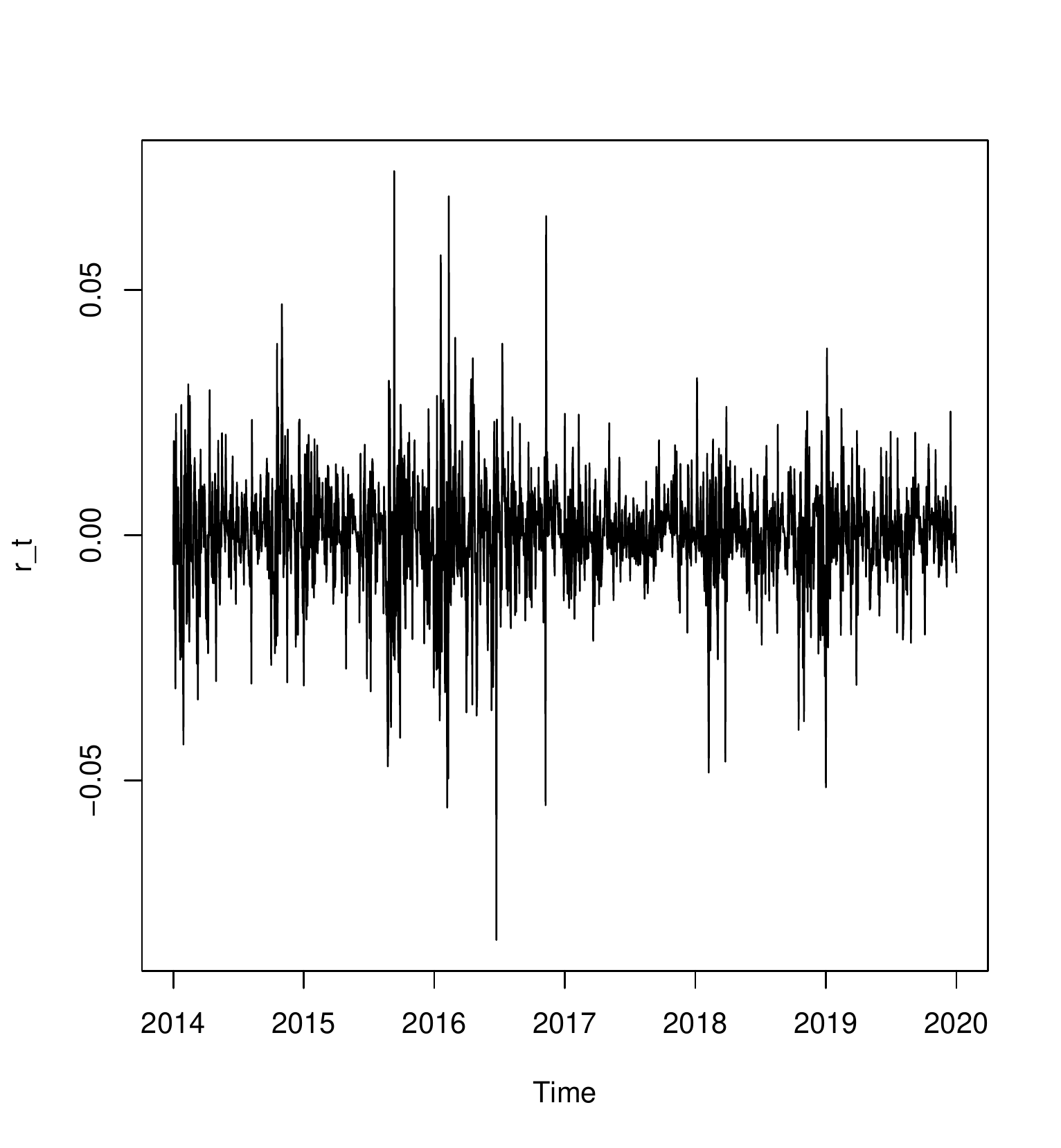}
            \caption{Time series of $r_t = \log (p_t/p_{t-1})$}\label{plot_lnikkei}
        \end{minipage}
\end{figure}

%\input Sec4.tex
%\section{Application :Analysis of Nikkei Stock Average}
In this section, we analyze the closing prices of the Nikkei Stock Average from January 1, 2014 to December 31, 2019, using tsQRF and compare the results with the WNW estimator. 
The data contained missing values, such as weekends and holidays, so we removed these missing values from the data. 
The data also show a long-term increasing trend (Figure \ref{plot_nikkei}), and we transform the price $p_t$ at time $t$ to $r_t = \log\frac{p_t}{p_{t-1}}$ (Figure \ref{plot_lnikkei}). 
In fact, the ADF test yields a p-value of $0.335$, so the null hypothesis of ``it is a unit root process" cannot be rejected. 
In this analysis, we used the first four years (length: 976) as training data and the remaining two years (length: 489) as test data.

Because the true conditional quantile values cannot be observed, the empirical coverage rate, which is defined as follows, is used to evaluate the accuracy of the quantile estimation and prediction: 

\begin{align*}
    \hat{\theta}&=\frac{N}{T}=\frac{1}{T}\sum_{t=1}^{T}\bm{1}_{\{Y_t\le \hat q^\tau(\bm x_t)\}}\quad \mbox{(training data)}, \\
    \tilde{\theta}&=\frac{N}{T'}=\frac{1}{T'}\sum_{t=1}^{T'}\bm{1}_{\{Y_{T+t}\le \hat q^\tau(\bm x_{T+t})\}}\quad \mbox{(test data)}.
\end{align*}

We fit the kernel quantile regression (WNW) and tsQRF with the orders $p=2$ and $p=5$ and we set the parameters to the same values as those used in Section 3 for functions \textbf{grf} and \textbf{np} in R packages.

Tables \ref{LR_nikkei_train} and \ref{LR_nikkei_test} show the estimated empirical coverage rates for the WNW and tsQRF for the training and test data, respectively. 
The values $\hat{\theta}$ and $\tilde{\theta}$, which are closer to $\tau$ between the WNW and tsQRF, are shown in red.

From Table \ref{LR_nikkei_train}, tsQRF provides more conservative result compared to WNW for the training data, and there is no significant difference in the accuracy of $\hat\theta$ and $\tilde\theta$ between $p=2$ and $p=5$. 
Table \ref{LR_nikkei_test} shows the prediction results for the test data, and tsQRF gives better estimated values for $\tau$ than WNW.

\begin{table}
\centering
\caption{$\hat{\theta}$ % and p-value of likelihood test
 on Training data}\label{LR_nikkei_train}
\begin{tabular}{cccccccccccc}
  \hline
  &&\multicolumn{5}{c}{$p=2$}&\multicolumn{5}{c}{$p=5$} \\
  \cmidrule(lr){3-7} \cmidrule(lr){8-12}
  &Model &2.5\% & 10\% & 50\% & 90\% & 97.5\%&2.5\% & 10\% & 50\% & 90\% & 97.5\%  \\  
  \hline
  $\hat{\theta}$& WNW &\textcolor{red}{0.027} & \textcolor{red}{0.102} & \textcolor{red}{0.501} & \textcolor{red}{0.891} & \textcolor{red}{0.969}&\textcolor{red}{0.024} & \textcolor{red}{0.097} & \textcolor{red}{0.496} & \textcolor{red}{0.894} & \textcolor{red}{0.974}  \\
   %& CAViaR & &&&&&&&&&  \\
  & tsQRF &0.000 & 0.041 & 0.463 & 0.880 & 0.969&0.000 & 0.037 & 0.462 & 0.890 & 0.970  \\
  \hline
\end{tabular}
\end{table}

\begin{table}
\centering
\caption{$\tilde{\theta}$ %and p-values of likelihood test 
on Test data}\label{LR_nikkei_test}
\begin{tabular}{cccccccccccc}
  \hline
  &&\multicolumn{5}{c}{$p=2$}&\multicolumn{5}{c}{$p=5$} \\
  \cmidrule(lr){3-7} \cmidrule(lr){8-12}
  &Model&2.5\% & 10\% & 50\% & 90\% & 97.5\%&2.5\% & 10\% & 50\% & 90\% & 97.5\%  \\  
  \hline
  $\tilde{\theta}$& WNW &0.012 & 0.084 & \textcolor{red}{0.497} & 0.926 & 0.988&0.016 & 0.086 & 0.487 & 0.926 & \textcolor{red}{0.982} \\ 
   %& CAViaR & &&&&&&&&&  \\
  & tsQRF &\textcolor{red}{0.022} & \textcolor{red}{0.108} & 0.493 & \textcolor{red}{0.922} & \textcolor{red}{0.984}&\textcolor{red}{0.018} & \textcolor{red}{0.086} & \textcolor{red}{0.497} & \textcolor{red}{0.924} & 0.988 \\ 
  \hline
%  p-value & kernel QR &\textcolor{blue}{0.046} & 0.222 & 0.892 & \textcolor{blue}{0.042} & \textcolor{blue}{0.046}&0.192 & 0.288 & 0.557 & \textcolor{blue}{0.042} & 0.327 \\ 
  %& CAViaR &&&&&&&&&&  \\
%  & tsQRF &0.718 & 0.541 & 0.752 & 0.089 & 0.192&0.327 & 0.288 & 0.892 & 0.062 & \textcolor{blue}{0.046} \\ 
%   \hline
\end{tabular}
\end{table}

\begin{figure}[htbp]
        \begin{minipage}[b]{0.46\textwidth}
            \centering
            \includegraphics[width=\columnwidth]{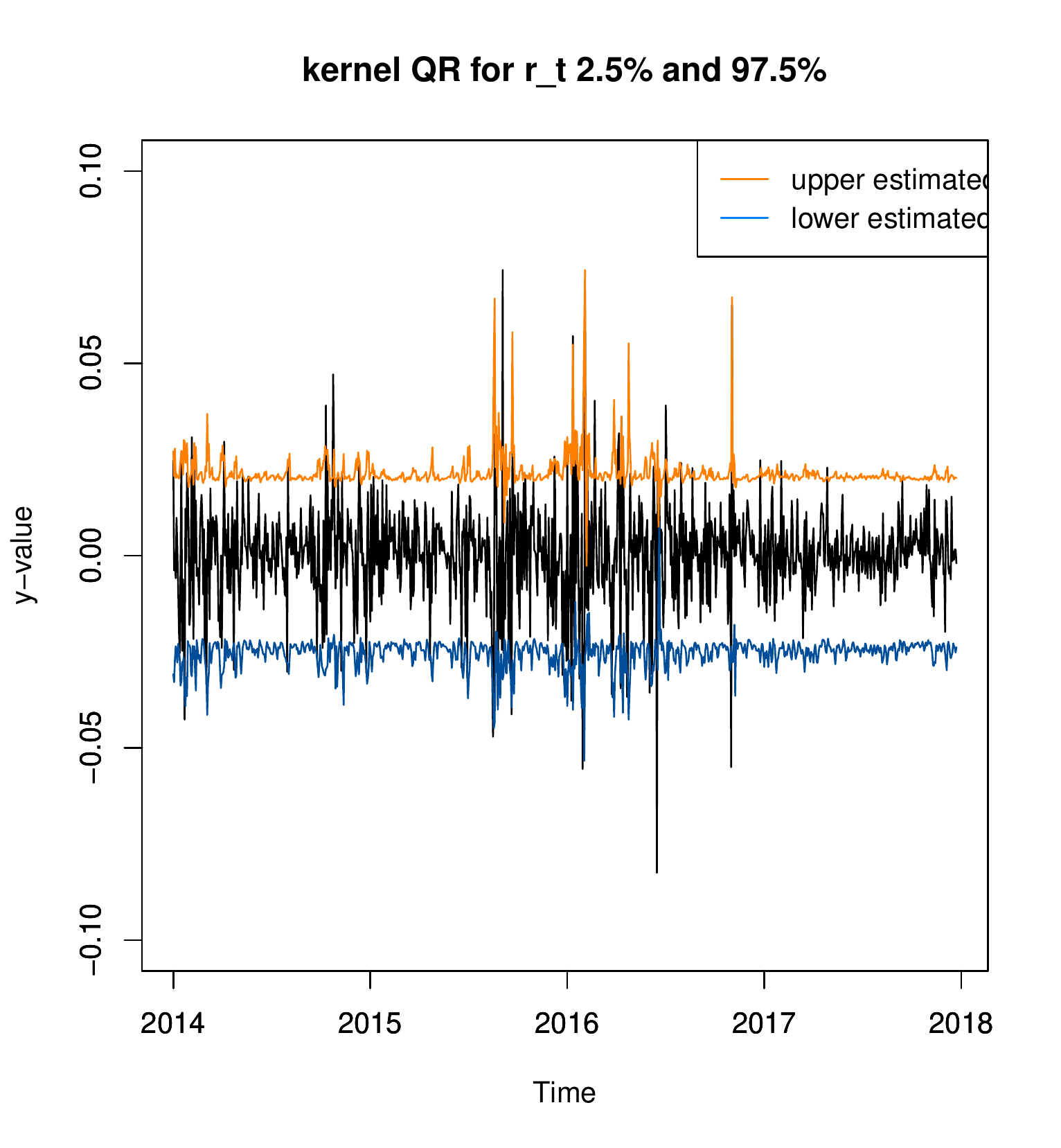}
            \caption{WNW ($p=5$)}\label{2.5kernel}
        \end{minipage}
        \hspace{0.04\textwidth} % ここで隙間作成
        \begin{minipage}[b]{0.46\textwidth}
            \centering
            \includegraphics[ width=\columnwidth]{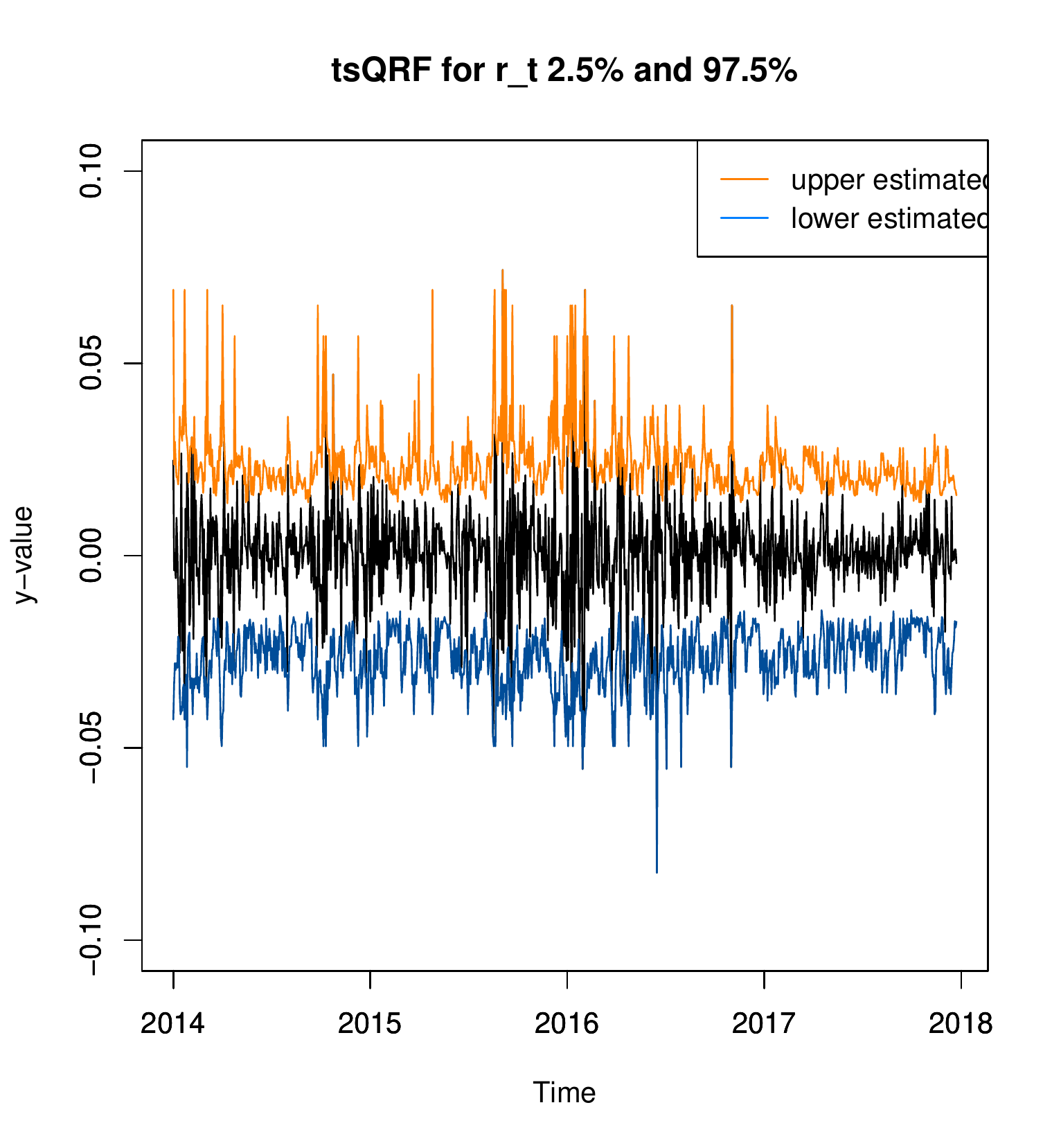}
            \caption{tsQRF ($p=5$)}\label{2.5grf}
        \end{minipage}
 %       \caption{$r_t$と分位点の変動（学習データ，$\tau=2.5\%,97.5\%$）}\label{2.5train}
\end{figure}

\begin{figure}[htbp]
        \begin{minipage}[b]{0.46\textwidth}
            \centering
            \includegraphics[width=\columnwidth]{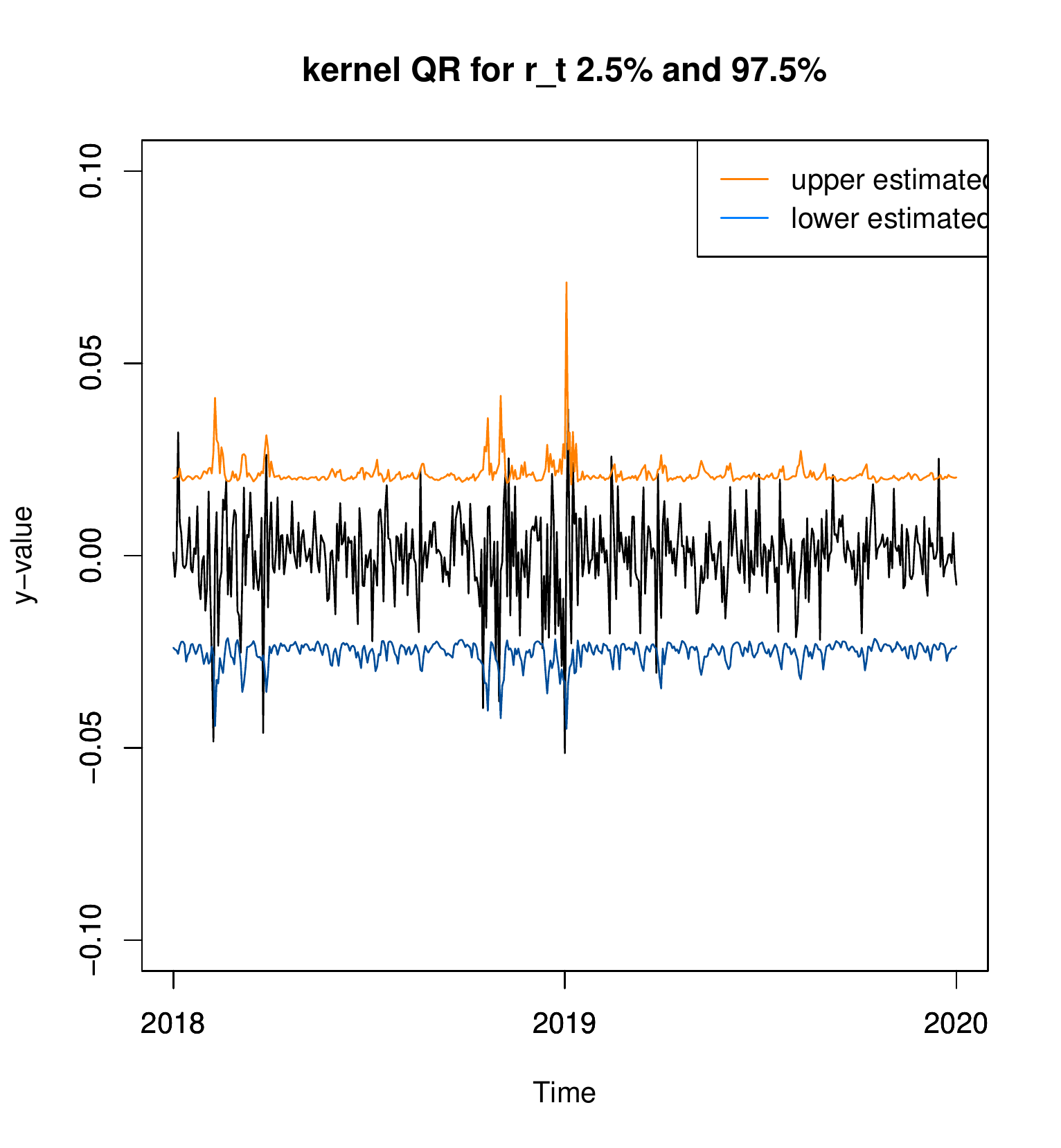}
            \caption{WNW ($p=5$)}\label{2.5kernel_test}
        \end{minipage}
        \hspace{0.04\textwidth} % ここで隙間作成
        \begin{minipage}[b]{0.46\textwidth}
            \centering
            \includegraphics[width=\columnwidth]{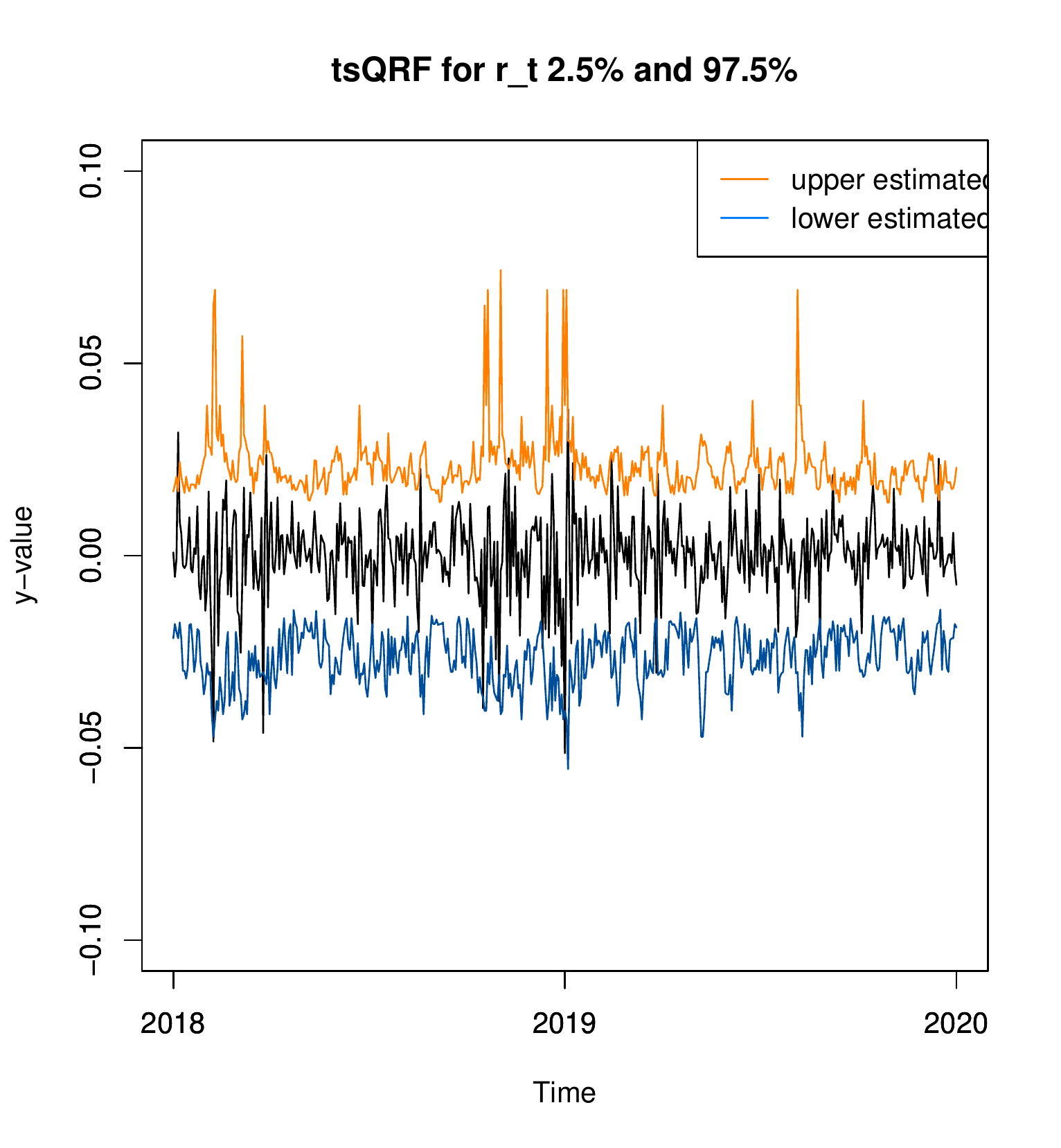}
            \caption{tsQRF ($p=5$)}\label{2.5grf_test}
        \end{minipage}
%        \caption{$r_t$と分位点の変動（テストデータ，$\tau=2.5\%,97.5\%$）}\label{2.5test}
\end{figure}

Figures \ref{2.5kernel} - \ref{2.5grf_test} illustrate the variation in the estimated quantiles of $r_t$ at $\tau=2.5\%,97.5\%$ for each estimation method, with $p=5$ for the training and test data.

Figures \ref{2.5kernel} and \ref{2.5grf} show that tsQRF is more sensitive to $r_t$ fluctuations than WNW. 
WNW does not capture the shocks of $r_t$. 
The same trend is also observed for the test data in Figures \ref{2.5kernel_test} and \ref{2.5grf_test}.

Figures \ref{2.5grf} and \ref{2.5grf_test} show that the quantile function estimated by tsQRF captures the variation of $r_t$. 
Even when there are large fluctuations in the time series, $r_t$ rarely exceeds the estimated 2.5\% and 97.5\% of points. 
Therefore, tsQRF is more sensitive to fluctuations in the time series than the WNW. 
For example, on September 9, 2015, concerns about the economy eased worldwide. 
Consequently, European and U.S. stocks rose the day before. 
The WNW could not capture the shock caused by this effect, whereas the tsQRF captured the effect properly. 
On June 24, 2016, supporters of leaving the European Union (EU) won the referendum in the United Kingdom. 
Consequently, there were concerns about the negative impact on the global economy, and the Nikkei Stock Average dropped sharply. 
Although tsQRF adequately captured these fluctuations, WNW did not. 

Based on these results, especially for data with large volatility fluctuations, tsQRF can capture the change in time series data with high sensitivity and may contribute to preventing underestimation of the risk rather than WNW.

\section{Summary and Future Work}\label{sec5}
We applied the Generalized Random Forests (GRF) proposed by Athey et al. (2019) %\cite{athey2019generalized}
 to quantile regression for time series data. 
Although theoretical confirmation has not been considered for their use in a time series setting, we derived the uniform consistency of the estimated function under mild conditions. 
Davis and Nielsen (2020) %\cite{davis2020modeling}
 also discussed the estimation problem using random forests (RF) for time series data, but the construction procedure of the RF treated by the GRF was essentially different, and different ideas were used throughout the theoretical proof. 
In addition, simulations and real data analyses were conducted. 
In the simulation, the accuracy of the conditional quantile estimation was evaluated under some time series models. 
In the real data using the Nikkei Stock Average, it was demonstrated that our estimator is more sensitive than the others in terms of  volatility and can prevent the underestimation of risk. 

However many challenges remain. 
If %Assumption 
\ref{A2} %(A-2)
 is relaxed, the range of applicable models can be expanded, including the traditional AR model. 
The model is expected to be extended to handle not only NLAR as in (\ref{eq2-1})%(2.1)
, but also ARCH-type models. 
Furthermore, in this study, the order $p$ is fixed; however, in practice, it should be determined using an information criterion or other methods. This is related to the variable selection problems in the GRF.
On the theoretical side, the discussion of asymptotic normality and asymptotic efficiency is the subject of future research. 
In particular, the efficiency involves the splitting procedure, and some methods, such as Neyman orthogonization, are expected to be effective.

%\color{black}
\paragraph{Acknowledgments}
This study was supported by JSPS KAKENHI Grant Number JP 21K11793. 
We would like to thank Editage (www.editage.com) for English language editing.

\bibliographystyle{plain}
%\bibliographystyle{plain.bst} %http://m0kichiazuma416.blog.fc2.com/blog-entry-44.html 
%\bibliography{biblio.bib}

%1. まず S-N-S.tex を pLaTex(ptex2pdf) でタイプセット（S-N-S.aux ができる）
%2. 出てきた S-N-S.aux を pBibTex で処理（S-N-S.bbl ができる）
%3. 文献リスト S-N-S.bbl ができたことを確認し，pLaTex(ptex2pdf) で S-N-S.tex を再度2回タイプセット

%数式番号のルール:  eq章-節-〇 (ex: eq4-2-1)

%\newpage
\appendix
\section{Proofs}\label{app}
\label{proofs}
Here, we present arguments leading up to our main result described in Section 2.
To prove Theorem \ref{thm1}, we first derive an upper bound of the moment of $\mathrm{diam}(\tilde L)$. 
Each leaf $\tilde L \in \tilde \Lambda$ can be expressed as $\tilde L := \bigotimes_{j=1}^p [r_j^-, r_j^+]$ based on a sequence $\{(r_j^-,r_j^+)\}_{j=1}^p$ with $0\le r_j^-< r_j^+ \le 1$ and $r_j^-, r_j^+ \in \mathcal{J}_s^Z$. 
Denoting the diameter of $\tilde L$ by $\mathrm{diam}(\tilde L)$, we can write
$$\mathrm{diam}(\tilde L) := \sup_{z', z'' \in \tilde L} \| z' -  z'' \| =\sqrt{\sum_{j=1}^p   \left| r_j^+ - r_j^-\right|^2}.$$
\begin{lemma}\label{lem1}
Under %Assumptions 
\ref{A1} -  \ref{A6}, we have
\begin{align*}
\mathbb{E}\left[ \left\| \mathrm{diam}(\tilde L) \right\|_{[0,1]^p}\right] = \mathbb{E}\left[ \sup_{z\in [0,1]^p }\left\{\mathrm{diam}(\tilde L(z)) \right\} \right] = O\left( s^{-\gamma} \right), 
\end{align*}
where 
$\displaystyle \gamma = \frac{\pi}{2} \frac{\log\left(\left(1-\omega \right)^{-1}\right)}{\log \left(\omega^{-1}\right)}$.
\end{lemma}

\paragraph{\textbf{Proof of Lemma \ref{lem1}}}
Let $c$ be the number of splits leading to any leaf $\tilde L\in \tilde \Lambda$, and let $c_j$ be the number of these splits along the $j$-th coordinate. 
Then, from %Assuptions 
\ref{A3} and \ref{A5}, 
\begin{align}
\mathrm{Binom}\left(c,p_j \right) \sim c_j \overset{d}{\ge} \underline c_j \sim \mathrm{Binom}\left(\underline{c},\pi \right), \label{eq4-1-1} 
\end{align}
where $\displaystyle \underline{c}= \frac{\log (s_*/(2k-1))}{\log \left(\omega^{-1}\right)}$ with 
$s_*= \left\lfloor \frac{s}{2} \right\rfloor$ (see, (31) of  Wager and Athey (2018)%\cite{wager2018estimation}
).  

By (\ref{eq4-1-1}) and Chernoff's inequality, for any $\{\delta_s \in(0,1)\} $ with $\delta_s \downarrow 0$, 
\begin{align}
\mathbb P\left[c_j \le (1+\delta_s)\frac{\underline c\pi}{2} \right]
&\le  \mathbb P\left[\underline c_j \le (1+\delta_s)\frac{\underline c\pi}{2} \right]
=  \mathbb P\left[\underline c_j \le \left(1 - \frac{1-\delta_s}{2}\right) \underline c\pi \right]\nonumber\\
&\le \exp\left( -\frac{1}{2}\left(\frac{1-\delta_s}{2}\right)^2 \underline c\pi \right)
= O\left( \exp\left(-\frac{\underline c \pi }{8}\right) \right).\label{eq4-1-2}
\end{align}
Since  
$$ \frac{\underline c\pi}{2} = \gamma\frac{\log (s_*/(2k-1))}{\log \left((1-\omega)^{-1}\right) }\quad \mathrm{and} \quad -\frac{ \underline c \pi }{8} < -\gamma \log (s_*/(2k-1))$$
from %Assumption
 \ref{A3}, (\ref{eq4-1-2}) yields that 
\begin{align}
\mathbb P\left[c_j \le \gamma(1+\delta_s) \frac{\log (s_*/(2k-1))}{\log \left((1-\omega)^{-1}\right) } \right]
=O\left(\exp\left(-\gamma \log \left(\frac{s_*}{2k-1}\right) \right) \right) = O\left (s^{-\gamma }\right). \label{eq4-1-3}
\end{align}

Let $\mathrm{diam}_j(\tilde L) = \mathrm{Leb}(\tilde L_j) =|r_j^+ - r_j^-|$ for any $\tilde L=\bigoplus_{j=1}^p [r_j^-,r_j^+] \in \tilde \Lambda$ with $\tilde L_j = [r_j^-, r_j^+]$. 
In what follows, %following the proof of Lemma 5.8 of Davis and Nielsen (2020), 
we will show that there exists some sequence $\{\delta_s \in (0,1)\}$ with $\delta_s \downarrow 0$ such that
\begin{align}
\mathrm{diam}_j(\tilde L) \le \left( 1-\omega\right)^{c_j/(1+\delta_s)}, \label{eq4-1-4}
\end{align}
with probability one. 
Let $\ell^{(j)} \in \{1\ldots, c_j\}$ be the depth of the tree $\tilde \Lambda$ splitted along the $j$-th coordinate, let $\tilde L^{\ell^{(j)}-}= \bigotimes_{i=1}^p \tilde L^{\ell^{(j)}-}_i$ and $\tilde L^{\ell^{(j)}}=\bigotimes_{i=1}^p \tilde L^{\ell^{(j)}}_i$ be the rectangles before and after the $\ell^{(j)}$-th splitting, respectively.   
Then, we can write
\begin{align}
\mathrm{diam}_j(\tilde L) 
%= \mathrm{diam}_j(\tilde L^{c_j}) 
=\mathrm{diam}_j(\tilde L^{1-}) \prod_{\ell^{(j)}=1}^{c_j} \frac{\mathrm{diam}_j(\tilde L^{\ell^{(j)}})}{\mathrm{diam}_j(\tilde L^{\ell^{(j)}-})}
=\prod_{\ell^{(j)}=1}^{c_j} \frac{\mathrm{diam}_j(\tilde L^{\ell^{(j)}})}{\mathrm{diam}_j(\tilde L^{\ell^{(j)}-})}.
\label{eq4-1-5}
\end{align}
For each rectangle $\tilde L^{\ell} = \bigotimes_{j=1}^p \tilde L_j^{\ell}$ with $ \tilde L_j^{\ell}= [r_j^{\ell -}, r_j^{\ell +}] \subset [0,1]$, there exists  $L^{\ell} = \bigotimes_{j=1}^p L_j^{\ell}$ with $L_j^{\ell} = [v_j^{\ell -}, v_j^{\ell +}] \subset \mathbb R$ such that $\tilde L_j^{\ell} = [F_h(v_j^{\ell -}), F_h(v_j^{\ell +})] =: F_h (L_j^{\ell})$ for $j=1,\ldots, p$, which implies that we have
\begin{align*}
\frac{\mathrm{diam}_j(\tilde L^{\ell^{(j)}})}{\mathrm{diam}_j(\tilde L^{\ell^{(j)}-})} 
= \frac{\mathrm{Leb}(F_h (L^{\ell^{(j)}}_j))}{\mathrm{Leb}(F_h (L^{\ell^{(j)}-}_j))} 
= \frac{\mathrm{Leb}(\iota_h (L^{\ell^{(j)}}))}{\mathrm{Leb}(\iota_h (L^{\ell^{(j)}-}))} 
= 1- \frac{\mathrm{Leb}(\iota_h (L^{\ell^{(j)}-}\setminus L^{\ell^{(j)}}))}{\mathrm{Leb}(\iota_h (L^{\ell^{(j)}-}))},
\end{align*}
where $\iota_h (L^{\ell}) = \bigotimes_{j=1}^p F_h(L_j^{\ell})$. 
In additon, from (5.45) and (5.46) of Davis and Nielsen (2020)%\cite{davis2020modeling}
, we have  
\begin{align}
0< \zeta^{-2}\omega \le \frac{\mathrm{diam}_j(\tilde L^{\ell^{(j)}})}{\mathrm{diam}_j(\tilde L^{\ell^{(j)}-})}
&\le 1- \zeta^{-2} \omega <1, \label{eq4-1-6}
\end{align}
where $\zeta= \bar{\zeta}^p$. 
Since $c_j\to \infty$ with probablity one, from (\ref{eq4-1-5}), (\ref{eq4-1-6}) and Glivenko-Cantelli theorem for ergodic process, it follows 
\begin{align*}
\frac{1}{c_j} \log \left(\mathrm{diam}_j(\tilde L)\right)  
= \frac{1}{c_j} \sum_{\ell(j)=1}^{c_j} \log \left(\frac{\mathrm{diam}_j(\tilde L^{\ell^{(j)}})}{\mathrm{diam}_j(\tilde L^{\ell^{(j)}-})} \right) %\nonumber \\
&\stackrel{a.s.}{\to} 
%\lim_{c_j\to \infty} 
\mathbb E \left[ \log \frac{\mathrm{diam}_j(\tilde L^{1})}{\mathrm{diam}_j(\tilde L^{1-})} \right] 
%= %\lim_{c_j\to \infty} 
%\mathbb E \left[ \log \frac{ \sharp \tilde L^{1}}{ \sharp \tilde L^{1-}} \right] 
\in [\log \omega, \log (1-\omega)],  %\label{eq4-1-7}
\end{align*}
which implies that there exits some sequence $\{\delta_s \in (0,1)\}$ with $\delta_s \downarrow 0$ satisfying (\ref{eq4-1-4}). 
Therefore, from (\ref{eq4-1-3}) and (\ref{eq4-1-4}), we have
\begin{align}
\mathbb P\left[ \mathrm{diam}_j(\tilde L) \ge (s_{*}/(2k-1))^{-\gamma} \right]
&\le \mathbb P\left[ (1-\omega)^{c_j/(1+\delta_s)} \ge (s_{*}/(2k-1))^{-\gamma} \right]\nonumber \\
&=\mathbb P\left[c_j \le \gamma(1+\delta_s) \frac{\log (s_{*}/(2k-1))}{\log \left((1-\omega)^{-1}\right) } \right]
=O\left(s^{-\gamma}\right) \label{eq4-1-7} 
\end{align}
Let $\max_{j\in\{1,\ldots, p\}}\mathrm{diam}_j(\tilde L)=:\mathrm{diam}_{j_*}(\tilde L)$. Then, it follows
\begin{align*}
\mathbb E\left[\mathrm{diam}(\tilde L)\right]
&\le %\sqrt{p} \mathbb E\left[\mathrm{diam}_{j_*}(\tilde L)\right]\\
%&= 
\sqrt{p} \left\{\mathbb E\left[\mathrm{diam}_{j_*}(\tilde L)\bm 1_{\{\mathrm{diam}_{j_*}(\tilde L) \ge (s_{*}/(2k-1))^{-\gamma}\}} \right]+\mathbb E\left[\mathrm{diam}_{j_*}(\tilde L)\bm 1_{\{\mathrm{diam}_{j_*}(\tilde L) < (s_{*}/(2k-1))^{-\gamma}\}} \right] \right\}\\
%&\le \sqrt{p} \left\{\mathbb P\left[\mathrm{diam}_{j_*}(\tilde L) \ge (s_{*}/(2k-1))^{-\gamma} \right]+ \mathbb E\left[(s_{*}/(2k-1))^{-\gamma} \bm 1_{\{\mathrm{diam}_{j_*}(\tilde L) < (s_{*}/(2k-1))^{-\gamma}\}} \right] \right\}\\
%&= O\left((s_{*}/(2k-1))^{-\gamma}\right) \\
&= O\left(s^{-\gamma}\right). 
\end{align*}
Since the above result does not depend on the test point $z=\iota_h(x)$, we obtain
\begin{align*}
\mathbb E\left[\| \mathrm{diam}(\tilde L) \|_{[0,1]^p}\right]
=\mathbb E\left[ \sup_{z\in [0,1]^p }  \mathrm{diam}(\tilde L(z)) \right]
= O\left(s^{-\gamma}\right). \qed
\end{align*}
\ \\

From Lemma \ref{lem1}, the upper bound of the moment of $\mathrm{diam}(L(x))$ for the original feature space $\mathcal X \subset \mathbb R^p$ can be derived. 
For any leaf $L\in \Lambda$, denoting the diameter of $L\cap \mathcal X \neq \emptyset$ by $\mathrm{diam}(L)$, we define 
$$\mathrm{diam}(L) := \sup_{x', x'' \in L\cap \mathcal X} \| x' -  x'' \|. $$
\begin{corollary}\label{cor1}
Under %Assumptions 
\ref{A1} - \ref{A6}, we have
\begin{align*}
\mathbb{E}\left[ \left\| \mathrm{diam}(L)\right\|_{\mathcal X}\right] = \mathbb{E}\left[ \sup_{x\in \mathcal X}  \mathrm{diam}(L(x)) \right] = O\left( s^{-\gamma} \right). 
\end{align*}
\end{corollary}

\paragraph{\textbf{Proof of Corollary \ref{cor1}}}
Since the map $\iota_h$ is one-to-one, from (\ref{eq4-1-7}), there exists some $C>0$ such that
\begin{align*}
\mathbb P\left[ \mathrm{diam}_j(L) \ge (s_*/(2k-1))^{-\gamma}\right]\le C (s_*/(2k-1))^{-\gamma},
\end{align*}
and since $L\cap \mathcal X \subset \mathcal X$ is a compact set, there exists some $C'>0$ such that
\begin{align*}
\mathbb E\left[ \mathrm{diam}_j(L) \bm 1_{\{\mathrm{diam}_j(L(x)) \ge (s_*/(2k-1))^{-\gamma} \}} \right] \le C' (s_*/(2k-1))^{-\gamma},
\end{align*}
for any $x\in \mathcal X$, which implies that the proof is completed in the same as that of Lemma \ref{lem1}. \qed

\ \\
\quad Denoting the conditional expectation of $\psi_{q(x)}(Y_t)$ given $X_t$ by $M_{q(x)}(X_t) = \mathbb E\left[ \psi_{q(x)}(Y_t)|X_t\right]$, we define
\begin{align*}
\bar{\Psi}_T(q)(x) := \sum_{t=1}^T  \alpha_t(x) M_{q(x)} (X_t),
\end{align*}
and for two parameter $q,q' \in \mathcal Q$ define
\begin{align*}
\delta(q,q') := \left\| \Psi_T(q) - \bar{\Psi}_T(q) - \left\{\Psi_T(q') - \bar{\Psi}_T(q')\right\} \right\|_{\mathcal X}. 
\end{align*}
Then, for the moments of $\delta(q,q')$, we can derive the followings. 
\begin{lemma}\label{lem2}
Under %Assumptions 
\ref{A1} - \ref{A6}, for any $q,q'\in \mathcal Q$, there exsits some $C>0$ such that 
\begin{align*}
\mathbb E\left[\delta(q,q')\right] = O\left( \frac{s}{T} \right), \quad 
\mathrm {Var} \left(\delta(q,q') \right)  \le C\frac{s\log s}{T} \left\| q-q' \right\|_{\mathcal X}.
\end{align*}
\end{lemma}
Note that this result corresponds to Lemma 8 of Athey et al. (2019) 
%\cite{athey2019generalized}
 in i.i.d. case. 
You can see that there exits a difference between i.i.d. case and dependent case. 
\paragraph{\textbf{Proof of Lemma \ref{lem2}}} 
For simplicity, we drop the index $x\in \mathcal X$ in  $q(x),\alpha_t(x),L(x)$ and so on.  
Denote
$$ \mathcal E_q (X_t,Y_t) := \psi_{q}(Y_t) - M_{q}(X_t), \quad \Delta \mathcal E_{q,q'}(X_t, Y_t) := \mathcal E_q (X_t,Y_t)- \mathcal E_{q'} (X_t,Y_t).$$
Then, we can write
\begin{align*}
\mathbb E\left[ \delta(q,q')\right] 
&= \frac{1}{|\mathcal{A}_s|}\sum_{A\in \mathcal{A}_s} \sum_{t\in A^{\mathcal I}} \mathbb E\left[ \frac{\bm 1_{\{ X_t\in L\}}}{\sharp L} \Delta \mathcal E_{q,q'} (X_t, Y_t) \right].
%&= \frac{1}{|\mathcal{A}_s|}\sum_{A\in \mathcal{A}_s} \sum_{t\in A^{\mathcal I}} \mathbb E\left[ \mathbb E \left[ \frac{\bm 1_{\{ X_t\in L\}}}{\sharp L} \Delta \mathcal E_{q,q'} (X_t, Y_t) |\mathcal I_s^X, \Lambda^{(A)}\right] | \Lambda^{(A)}\right].
\end{align*}
Fix $A=\{A^{\mathcal I}, A^{\mathcal J}\}\in \mathcal A_s$; $L\in \Lambda^{(A)}$ where $\Lambda^{(A)}$ is a partitioning $\Lambda$ based on $A\in \mathcal A_s$; and $t\in A^{\mathcal I}$.  
If $|t - t^{\mathcal J}|> p+1$ is satisfied for all $t^{\mathcal J} \in A^{\mathcal J}$, it follows from the $p$-dependency of $\{Y_t\}$ 
\begin{align}
\mathbb E \left[ \frac{\bm 1_{\{ X_t\in L\}}}{\sharp L} \Delta \mathcal E_{q,q'} (X_t, Y_t) |\mathcal I_s^X, L\right] 
=\frac{\bm 1_{\{ X_t\in L\}}}{\sharp L} \left\{ \mathbb E \left[  \mathcal E_{q} (X_t, Y_t) |X_t\right] -  \mathbb E \left[\mathcal E_{q'} (X_t, Y_t) |X_t\right]\right\}
=0. \label{eqA-2-1}
\end{align}
Otherwise from %Assumption 
\ref{A5}, $|\Delta \mathcal E_{q,q'}(X_t, Y_t)|\le 2$, and the stationarity of $\{Y_t\}$, 
\begin{align}
\mathbb E \left[ \frac{\bm 1_{\{ X_t\in L\}}}{\sharp L} \Delta \mathcal E_{q,q'} (X_t, Y_t) |L\right] 
\le \frac{2}{k} \mathbb E \left[ \bm 1_{\{ X_t\in L\}} |L \right]
=\frac{2}{ks_*} \mathbb E \left[ \sum_{t\in A^{\mathcal I} }\bm 1_{\{ X_t\in L\}} |L \right] 
\le \frac{2(2k-1)}{ks_*}
\label{eqA-2-2}
\end{align}
where $s_*=\left\lfloor \frac{s}{2} \right\rfloor =|A^{\mathcal I}|$. 
Let $\mathcal{K}_{\ell}^{(1)}$ ($\ell =\pm 1,\ldots, \pm (T-1)$) be a multiset of $t\in \{1,\ldots, T\}$ satisfying that there exists an $A=\{A^{\mathcal I},A^{\mathcal J}\}\in \mathcal A_s$ such that $t \in A^{\mathcal I}$ and $t+\ell \in A^{\mathcal J}$. 
Note that if there exist $A_1=\{A_1^{\mathcal I},A_1^{\mathcal J}\}\in \mathcal A_s$ and $A_2=\{A_2^{\mathcal I},A_2^{\mathcal J}\}\in \mathcal A_s$ with $A_1\neq A_2$, $A_1^{\mathcal I}\ni t \in A_2^{\mathcal I}$ and  $A_1^{\mathcal J}\ni t+\ell \in A_2^{\mathcal J}$, the number of element $t$ in $\mathcal{K}_{\ell}^{(1)}$ is two (see an example in Table \ref{example_of_K}). 

\begin{table}
\centering
\caption{Example of $\mathcal{A}_{s}$ and $\mathcal{K}_{\ell}^{(1)}$ in case of $T=5, s=4$} \label{example_of_K}
\small%\tiny
\begin{tabular}{cccccccc}
  \hline
 $\mathcal{A}_s$ & $\{\{1,2\},\{3,4\}\}$ & $\{\{1,2\},\{3,5\}\}$ &$\{\{1,2\},\{4,5\}\}$ 
                        & $\{\{1,3\},\{2,4\}\}$ & $\{\{1,3\},\{2,5\}\}$  &$\{\{1,3\},\{4,5\}\}$ \\  
                        & $\{\{1,4\},\{2,3\}\}$ & $\{\{1,4\},\{2,5\}\}$ &$\{\{1,4\},\{3,5\}\}$ 
                        & $\{\{1,5\},\{2,3\}\}$ & $\{\{1,5\},\{2,4\}\}$  &$\{\{1,5\},\{3,4\}\}$ \\  
                        & $\{\{2,3\},\{1,4\}\}$ & $\{\{2,3\},\{1,5\}\}$ &$\{\{2,3\},\{4,5\}\}$ 
                        & $\{\{2,4\},\{1,3\}\}$ & $\{\{2,4\},\{1,5\}\}$  &$\{\{2,4\},\{3,5\}\}$ \\  
                        & $\{\{2,5\},\{1,3\}\}$ & $\{\{2,5\},\{1,4\}\}$ &$\{\{2,5\},\{3,4\}\}$ 
                        & $\{\{3,4\},\{1,2\}\}$ & $\{\{3,4\},\{1,5\}\}$  &$\{\{3,4\},\{2,5\}\}$ \\  
                        & $\{\{3,5\},\{1,2\}\}$ & $\{\{3,5\},\{1,4\}\}$ &$\{\{3,5\},\{2,4\}\}$ 
                        & $\{\{4,5\},\{1,2\}\}$ & $\{\{4,5\},\{1,3\}\}$  &$\{\{4,5\},\{2,3\}\}$ \\  
  \hline
 $\mathcal{K}_{1}^{(1)}$ & \multicolumn{6}{l}{ $\{1,1,1,1,1,1,2,2,2,2,2,2,3,3,3,3,3,3,4,4,4,4,4,4\}$ } \\
 $\mathcal{K}_{2}^{(1)}$ & \multicolumn{6}{l}{ $\{1,1,1,1,1,1,2,2,2,2,2,2,3,3,3,3,3,3\}$ } \\
 $\mathcal{K}_{3}^{(1)}$ & \multicolumn{6}{l}{ $\{1,1,1,1,1,1,2,2,2,2,2,2\}$ } \\
 $\mathcal{K}_{4}^{(1)}$ & \multicolumn{6}{l}{ $\{1,1,1,1,1,1\}$ } \\
 $\mathcal{K}_{-1}^{(1)}$ & \multicolumn{6}{l}{ $\{2,2,2,2,2,2,3,3,3,3,3,3,4,4,4,4,4,4,5,5,5,5,5,5\}$ } \\
 $\mathcal{K}_{-2}^{(1)}$ & \multicolumn{6}{l}{ $\{3,3,3,3,3,3,4,4,4,4,4,4,5,5,5,5,5,5\}$ } \\
 $\mathcal{K}_{-3}^{(1)}$ & \multicolumn{6}{l}{ $\{4,4,4,4,4,4,5,5,5,5,5,5\}$ } \\
 $\mathcal{K}_{-4}^{(1)}$ & \multicolumn{6}{l}{ $\{5,5,5,5,5,5\}$ } \\
   \hline
\end{tabular}
\end{table}   
Then, we can interpret
\begin{align*}
\sum_{A \in \mathcal A_s}  \sum_{t \in A^{\mathcal I}}  
\equiv  \sum_{|\ell|=1}^{T-1} \sum_{t \in \mathcal K_{\ell}^{(1)}} \left(
  \begin{array}{c}
    T-2  \\
    s-2 \\
  \end{array}
\right)\left(
  \begin{array}{c}
    s-2  \\
    s_*-1 \\
  \end{array}
\right)\left(
  \begin{array}{c}
    s-s_*-1  \\
    s-s_*-1 \\
  \end{array}
\right),
\end{align*}
where $\left(
  \begin{array}{c}
    T-2  \\
    s-2 \\
  \end{array}
\right)$ is the number of cases of all $A\in \mathcal A_s$ except for $\{t,t+\ell\}\in A$ and $\left(
  \begin{array}{c}
    s-2  \\
    s_*-1 \\
  \end{array}
\right)\left(
  \begin{array}{c}
    s-s_*-1  \\
    s-s_*-1 \\
  \end{array}
\right)$ is the number of cases of the divison for $s-2$ elements into $A^{\mathcal I}$ and $A^{\mathcal J}$ except for $t\in A^{\mathcal I}$ and $t+\ell \in A^{\mathcal J}$. 
Therefore, from (\ref{eqA-2-1}) and (\ref{eqA-2-2}), we have
\begin{align*}
\mathbb E\left[ \delta(q,q')\right] 
&\le \frac{1}{|\mathcal A_s|} \sum_{A \in \mathcal A_s} \sum_{t \in A^{\mathcal I}} \bm 1_{\{ \exists (t+ \ell) \in A^{\mathcal J},\ |\ell|\le p \}} \frac{2(2k-1)}{ks_*} \\
&= \frac{1}{|\mathcal A_s|} \sum_{|\ell|=1}^p \sum_{t \in \mathcal K_{\ell}^{(1)}} \left(
  \begin{array}{c}
    T-2  \\
    s-2 \\
  \end{array}
\right)\left(
  \begin{array}{c}
    s-2  \\
    s_*-1 \\
  \end{array}
\right)\frac{2(2k-1)}{ks_*}\\
&\le \left(
  \begin{array}{c}
    T  \\
    s \\
  \end{array}
\right)^{-1}\left(
  \begin{array}{c}
    s  \\
    s_* \\
  \end{array}
\right)^{-1} (2pT) \left(
  \begin{array}{c}
    T-2  \\
    s-2 \\
  \end{array}
\right)\left(
  \begin{array}{c}
    s-2  \\
    s_*-1 \\
  \end{array}
\right)\frac{2(2k-1)}{ks_*}\\
&=\frac{s(s-1)}{T(T-1)} \frac{s_*(s-s_*)}{s(s-1)} (2pT) \frac{2(2k-1)}{ks_*}\\
&= O\left( \frac{s}{T} \right),
%\label{eqA-2-3}
\end{align*}

\noindent which is completed the proof of the first claim. 

To prove the second claim, we denote the subsamples $(\mathcal I_s^i, \mathcal J_s^i)$ corresponds to $A_i=\{A_i^{\mathcal I}, A_i^{\mathcal J}\}\in \mathcal A_s$ for $i=1,2$, and introduce  
\begin{align*}
\Delta \mathcal T_{q,q'}(\mathcal I_s^i, \mathcal J_s^i) 
:= \sum_{t\in A_i^{\mathcal I}}\alpha_{A_i,t} \Delta \mathcal E_{q,q'}(X_t, Y_t) .
\end{align*}
Let the conditional covariance of $\Delta \mathcal T_{q,q'}$ for $A_1$ $and$ $A_2$ given $L_1\in \Lambda^{(A_1)}$ and  $L_2\in \Lambda^{(A_2)}$  by  $\mathrm{Cov}_{\Lambda}$. Then, we can write
\begin{align*}
&\mathrm {Cov}_{\Lambda} \left( \Delta \mathcal T_{q,q'}(\mathcal I_s^1, \mathcal J_s^1), \Delta \mathcal T_{q,q'}(\mathcal I_s^2, \mathcal J_s^2)\right) \nonumber \\
&=  \sum_{(t_1,t_2) \in A_1^{\mathcal I}\times A_2^{\mathcal I}} \mathrm{Cov}_{\Lambda} \left( \alpha_{A_1,t_1} \Delta \mathcal E_{q,q'} (X_{t_1}, Y_{t_1}), \alpha_{A_2,t_2} \Delta \mathcal E_{q,q'} (X_{t_2}, Y_{t_2})\right). %\label{eq4-2-2}
\end{align*}
For each $(t_1,t_2) \in A_1^{\mathcal I}\times A_2^{\mathcal I}$, from the continuously of $F_{\epsilon}$ and %Assumption 
\ref{A5}, there exists some $\tilde C\ge 0$ such that
\begin{align*}
\mathrm{Cov}_{\Lambda} \left( \alpha_{A_1,t_1} \Delta \mathcal E_{q,q'} (X_{t_1}, Y_{t_1}), \alpha_{A_2,t_2} \Delta \mathcal E_{q,q'} (X_{t_2}, Y_{t_2})\right)
%&\le \frac{1}{k^2}\left|\mathrm{Cov}_{\Lambda} \left( \bm 1_{\{X_{t_1} \in L_1\}} \Delta \mathcal E_{q,q'} (X_{t_1}, Y_{t_1}), \bm 1_{\{X_{t_2}\in L_2\}} \Delta \mathcal E_{q,q'} (X_{t_2}, Y_{t_2})\right)\right|\\
&\le %\frac{\tilde C}{k^2}\left| \mathrm{Cov}_{\Lambda} \left( \bm 1_{\{X_{t_1} \in L_1\}} , \bm 1_{\{X_{t_2}\in L_2\}} \right) \right| |q-q'|\\
%&= 
\frac{\tilde C}{k^2}\left| R_{\Lambda}(t_1-t_2) \right| \|q-q'\|_{\mathcal X},
\end{align*}
where $R_{\Lambda}(t_1-t_2) := \mathrm{Cov}_{\Lambda} \left( \bm 1_{\{X_{t_1} \in L_1\}} , \bm 1_{\{X_{t_2}\in L_2\}} \right)$.
Hence, we have
\begin{align}
%\mathrm {Var} \left( \delta(q,q')\right) 
%&= 
&\frac{1}{|\mathcal A_s|^2} \sum_{A_1, A_2 \in \mathcal A_s} \mathrm{Cov}_{\Lambda}\left( \Delta \mathcal T_{q,q'}(\mathcal I_s^1, \mathcal J_s^1), \Delta \mathcal T_{q,q'}(\mathcal I_s^2, \mathcal J_s^2)\right)\nonumber \\
&\le \left( \frac{\tilde C}{k^2} \|q-q'\|_{\mathcal X}\right)  \frac{1}{|\mathcal A_s|^2} \sum_{A_1, A_2 \in \mathcal A_s}  \sum_{(t_1,t_2) \in A_1^{\mathcal I}\times A_2^{\mathcal I}} \left| R_{\Lambda}(t_1-t_2) \right|. \label{eq4-2-3}
\end{align}

Let $\mathcal{K}_{\ell}^{(2)}(\ell =0,\pm 1,\ldots, \pm (T-1))$ be a multiset of $t\in \{1,\ldots,T\}$ satisfying that there exists an $(A_1,A_2)=(\{A_1^{\mathcal I}, A_1^{\mathcal J}\}, \{A_2^{\mathcal I}, A_2^{\mathcal J}\})\in \mathcal A_s^2$ such that $t \in A_1^{\mathcal I}$ and $t+\ell \in A_2^{\mathcal I}$. 
Then, we can interpret 
\begin{align*}
\sum_{A_1,A_2 \in \mathcal A_s}  \sum_{(t_1,t_2) \in A_1^{\mathcal I}\times A_2^{\mathcal I}}  
\equiv  \sum_{|\ell|=0}^{T-1} \sum_{t \in \mathcal K_{\ell}^{(2)}} \left(
  \begin{array}{c}
    T-1  \\
    s-1 \\
  \end{array}
\right)^2 \left(
  \begin{array}{c}
    s-1  \\
    s_*-1 \\
  \end{array}
\right)^2\left(
  \begin{array}{c}
    s-s_*  \\
    s-s_* \\
  \end{array}
\right)^2,
\end{align*}
where $\left(
  \begin{array}{c}
    T-1  \\
    s-1 \\
  \end{array}
\right)^2$ is the number of cases of all $A_1,A_2\in \mathcal A_s$ except for $t\in A_1, t+\ell \in A_2$ and $\left(
  \begin{array}{c}
    s-1  \\
    s_*-1 \\
  \end{array}
\right)\left(
  \begin{array}{c}
    s-s_*  \\
    s-s_* \\
  \end{array}
\right)^2$ is the number of cases of the divison for $s-1$ elements into $A_1^{\mathcal I} (A_2^{\mathcal I})$ and $A_1^{\mathcal J}(A_2^{\mathcal J})$ except for $t\in A_1^{\mathcal I}(t+\ell \in A_2^{\mathcal I})$. 
Therefore, we have
\begin{align}
&\frac{1}{|\mathcal A_s|^2} \sum_{A_1, A_2 \in \mathcal A_s}  \sum_{(t_1,t_2) \in A_1^{\mathcal I}\times A_2^{\mathcal I}} \left| R_{\Lambda}(t_1-t_2) \right| \nonumber \\
&=\frac{1}{|\mathcal A_s|^2} \sum_{|\ell| = 0}^{T-1} \sum_{t  \in \mathcal K_{\ell}^{(2)}} \left(
  \begin{array}{c}
    T-1  \\
    s-1 \\
  \end{array}
\right)^{2}\left(
  \begin{array}{c}
    s-1  \\
    s_*-1 \\
  \end{array}
\right)^{2}\left| R_{\Lambda}(\ell) \right| \nonumber \\
&\le \left(
  \begin{array}{c}
    T  \\
    s \\
  \end{array}
\right)^{-2} \left(
  \begin{array}{c}
    s  \\
    s_* \\
  \end{array}
\right)^{-2} T
\left(
  \begin{array}{c}
    T-1  \\
    s-1 \\
  \end{array}
\right)^{2}
\left(
  \begin{array}{c}
   s-1  \\
    s_*-1 \\
  \end{array}
\right)^{2}  \sum_{|\ell| = 0}^{T-1}  \left| R_{\Lambda} (\ell) \right| \nonumber \\
&= \frac{s^2}{T^2} \frac{s_*^2}{s^2} T \sum_{|\ell| = 0}^{T-1}  \left| R_{\Lambda} (\ell) \right| \nonumber \\
&= \frac{s_*^2}{T} \sum_{|\ell| = 0}^{T-1}  \left| R_{\Lambda} (\ell) \right| . \label{eq4-2-4}
\end{align}

On the other hand, as we state in Remark 1, $\{Y_t\}_{t\in \mathbb N}$ is exponentially $\alpha$-mixing and %$X_t=(Y_{t-1},\ldots, Y_{t-p})$の関数で定義される
$\{\bm 1_{\{X_t \in L\}}\}_{t\in \mathbb N}$ is also exponentially $\alpha$-mixing. 
Therefore, denoting $\mu(L)= \mathbb P\left[ X_t \in L | L \right]$, we have
$|R_{\Lambda}(\ell)| = \left| \mathrm{Cov}_{\Lambda} \left( \bm 1_{\{X_t \in L\}}, \bm 1_{\{ X_{t+\ell} \in L\}}\right) \right| \le \min\left\{ \alpha (\ell), \mu (L) \right\}$ 
with the mixing coefficient $\alpha(\ell)\lesssim e^{-\ell}$. 
Since 
\begin{align*}
\mu(L):= \mathbb P\left[ X_t \in L | L \right]=\mathbb E\left[ \bm 1_{\{X_t \in L\}}|L \right] = \frac{1}{s_*}  \mathbb E\left[ \sum_{t \in A^{\mathcal I}} \bm 1_{\{ X_t \in L\}}|L\right] \le \frac{2k-1}{s_*} 
\end{align*}  
from %Assumption 
\ref{A5} and the stationarity of $\{X_t\}$, we have 
\begin{align}
\sum_{\ell=0}^{\infty} |R_{\Lambda}(\ell)| 
&\lesssim \sum_{\ell \le \lfloor \log \frac{s_*}{2k-1} \rfloor } \frac{2k-1}{s_*} + \sum_{\ell > \lfloor \log \frac{s_*}{2k-1} \rfloor} e^{-\ell} %\nonumber \\
\lesssim  \frac{\log \frac{s_*}{2k-1} }{s_*/(2k-1)} +  e^{-\log \frac{s_*}{2k-1}} \frac{1}{1-e}\nonumber\\
&=O\left(\frac{\log s}{s} \right). \label{eq4-2-6}
\end{align}  

\noindent The above evaluation does not depend on the condition $\Lambda^{(A_1)}$ and  $\Lambda^{(A_2)}$, and the dependency of $\Lambda^{(A_1)}$ and  $\Lambda^{(A_2)}$ is negligible with high probability from the first claim, so that 
\begin{align*}
\mathrm {Var} \left( \delta(q,q')\right) 
\lesssim \frac{1}{|\mathcal A_s|^2} \sum_{A_1, A_2 \in \mathcal A_s} \mathrm{Cov}_{\Lambda}\left( \Delta \mathcal T_{q,q'}(\mathcal I_s^1, \mathcal J_s^1), \Delta \mathcal T_{q,q'}(\mathcal I_s^2, \mathcal J_s^2)\right)
\end{align*}
which implies from (\ref{eq4-2-3})-(\ref{eq4-2-6}), that there exists som $C> 0$ such that
\begin{align*}
\mathrm {Var} \left( \delta(q,q')\right) 
\lesssim \left( \frac{\tilde C}{k^2} \|q-q'\|_{\mathcal X}\right)  \frac{s_*^2}{T} \left\{ R_{\Lambda}(0) + 2\sum_{\ell=1}^{\infty}|R_{\Lambda}(\ell)| \right\} 
\le C \frac{s\log s}{T} \|q-q'\|_{\mathcal X}.
\end{align*}
Although the above inequality is obtained under any fixed $x \in \mathcal X$, the test point $x$ depends only on $|q(x)-q'(x)|$ in the right hand side, and replacing this part by $\|q-q'\|_{\mathcal X}$, the proof of the second claim is compleated. \qed

\ \\
By Lemmas \ref{lem1} and \ref{lem2}, we can derive the uniform consistency of $\Psi_T$ for $\mathcal Q$. 
\begin{lemma}\label{lem3}
Under %Assumptions 
\ref{A1} -  \ref{A6}, we have
\begin{align*}
\left\| \Psi_T- \Psi \right\|_{\mathcal Q}\stackrel{p}{\to} 0 \quad \mathrm{as}\ T\to\infty
\end{align*}
\end{lemma}
\paragraph{\textbf{Proof of Lemma \ref{lem3}}}
From the triangle inequality, 
$$ \| \Psi_T - \Psi \|_{\mathcal Q} \le  \| \Psi_T - \bar{\Psi}_T \|_{\mathcal Q}+ \| \bar{\Psi}_T - \Psi \|_{\mathcal Q}=: I_1+I_2, $$
it is sufficient to show $I_j \stackrel{p}{\to} 0$ ($j=1,2$). 

For $I_1$, we consider the empirical process technique for U-statistics. 
Let any subsample $(\mathcal I_s, \mathcal J_s)=(\mathcal D_{A^{\mathcal I}},\mathcal D_{A^{\mathcal J}})$ for $A=\{A^{\mathcal I},A^{\mathcal J}\} \in \mathcal A_s$ be an random element taking value on the sample space $((\mathcal X\times \mathcal Y)^s, (\mathcal B_X\times \mathcal B_Y)^s)$ where $\mathcal B_X$ and $\mathcal B_Y$ are the Borel sets on $\mathbb R^p$ and $\mathbb R$, respectively. 
Define the empirical measure by $\mathbb P_{s,T}:= \frac{1}{|\mathcal A_s|} \sum_{A\in \mathcal A_s} \delta_{(\mathcal D_A^{\mathcal I},\mathcal D_A^{\mathcal J})}$
where $\delta$ is a Dirac measure. 
Then, for any $(\mathcal I_s,\mathcal J_s)$-mesurable function $f$, we can write
$$  \mathbb P_{s,T} f =  \frac{1}{|\mathcal A_s|} \sum_{A\in \mathcal A_s} f\left(\mathcal I_s,\mathcal J_s\right).$$
We also define $P$ as a probabilty measure for any $(\mathcal I_s,\mathcal J_s)$ with $Pf = \int f dP = \mathbb E\left[ f(\mathcal I_s,\mathcal J_s)\right].$
Let $y=q(x) \in \mathcal Y \subset \mathbb R$ be a parameter. 
We introduce an $(\mathcal I_s,\mathcal J_s)$-mesurable function $\mathcal T_y$ by 
$$ \mathcal T_{y} (\mathcal I_s,\mathcal J_s)(x) := \sum_{t\in A^{\mathcal I}} \alpha_{A,t}(x) \mathcal E_{y} (X_t, Y_t),$$
for any $x\in \mathcal X$, and denote $\mathcal F$ by the set of function $\mathcal T_{y}$. 
Moreover, for any signed meaure $G$, we difine 
$$ \|G\|_{\mathcal F} = \sup_{\mathcal T_y \in \mathcal F} \|G\mathcal T_y\|_{\mathcal X} = \sup_{y \in \mathcal Y} \|G\mathcal T_y\|_{\mathcal X}=\sup_{q \in \mathcal Q} \|G\mathcal T_{q(\cdot)}\|_{\mathcal X},$$
then we can write
\begin{align*}
\left\| \mathbb P_{s,T}- P \right\|_{\mathcal F}  
&= \sup_{q \in \mathcal Q} \left\| \Psi _T (q) - \bar{\Psi}_T(q) - \mathbb E \left[ \Psi _T (q) - \bar{\Psi}_T(q) \right] \right\|_{\mathcal X}.
\end{align*}
Under a fixed $\mathcal D_T$, for any $\epsilon>0$, we construct an $\epsilon$-bracket satisfying that the $L_1(\mathbb P_{s,T})$-norm is less than $\epsilon$. 
Define 
\begin{align*}
\mathcal T_{y}(\mathcal I_s, \mathcal J_s)(x)  &= \left\{ \sum_{t\in A^{\mathcal I}} \alpha_{A,t}(x)F_{\varepsilon}(y- g(X_t)) \right\} - \left\{ \sum_{t\in A^{\mathcal I}}\alpha_{A,t}(x) \bm 1_{\{Y_t \le y\}} \right\} \\
&=: \mathcal T_{y}^1(\mathcal I_s, \mathcal J_s)(x) - \mathcal T_{y}^2(\mathcal I_s, \mathcal J_s)(x)
\end{align*}
and $\epsilon_j =  j (\epsilon/2) $ for $j =0,\ldots, \lfloor 2/\epsilon  \rfloor$. 
For each $j \in \{ 0,\ldots, \lfloor 2/\epsilon  \rfloor\},\ i\in \{1,2\}$, we also define
$$ y_j^i := \inf \left\{ y\in \mathcal Y : \mathbb P_{s,T} \mathcal T_{y}^i(x) \ge \epsilon_{j},\  \forall x\in \mathcal X \right\}.$$ 

\noindent Since $\mathcal T_{y}^1$ and $ \mathcal T_{y}^2$ are monotonically non-decreasing functions with respect to $x$ under any fixed $\mathcal I_s $ and $ \mathcal J_s$, it follows $y_j^i\le y_{j'}^i$ if $j<j'$. 
Let $\{\tilde{y}_j\}_{j=0,\ldots, 2\lfloor 2/\epsilon  \rfloor} \equiv \{y_j^1\}\oplus \{y_j^2\}$ satisfying $\tilde y_j\le \tilde y_{j'}$ if $j<j'$ and let $\mathcal F_{j} =  \left\{ \mathcal T_y \in \mathcal F :  \tilde y_j \le y \le \tilde y_{j+1}^- \right\}$,  
where $\tilde y_{j+1}^-=\inf \{ y\in \mathcal Y : y \ge \tilde y_{j+1}\}$ and $\mathcal F_{j}=\emptyset$ if $\tilde y_j = \tilde y_{j+1}$. 
Then, for any $\mathcal T_y, \mathcal T_{y'} \in \mathcal F_{j}\setminus \emptyset$, and any $x\in \mathcal X$, there exists some $y_{j(i)}^i, y_{j(i)+1}^i \in \{y_j^i\}$, such that $y_{j(i)}^i \le \tilde y_j <\tilde y_{j+1}\le y_{j(i)+1}^i$, and 
\begin{align*}
\mathbb P_{s,T} |\mathcal T_{y}(x)-\mathcal T_{y'}(x) |  
&\le \sum_{i=1}^2 \mathbb P_{s,T} |\mathcal T_{y}^i(x)-\mathcal T_{y'}^i(x) | 
\le   \sum_{i=1}^2 \mathbb P_{s,T} \left(\mathcal T_{\tilde y_{j+1}^-}^i(x)-\mathcal T_{\tilde y_{j}}^i(x) \right)\\
&\le   \sum_{i=1}^2 \mathbb P_{s,T} \left(\mathcal T_{y_{j(i)+1}^-}^i(x)-\mathcal T_{y_{j(i)}}^i(x) \right)
<  \sum_{i=1}^2 \left( \epsilon_{j(i)+1} - \epsilon_{y_{j(i)}}\right)
 = \epsilon,
\end{align*}
which implies that $\mathcal F_j$ is $\epsilon$-bracket with respect to $L_1(\mathbb P_{s,T})$ norm, and the bracketing number which is the minimum number of brackets of radius at most $\epsilon$ required to cover $\mathcal F$ in terms of the $L_1(\mathbb P_{s,T})$ norm, satisfies
\begin{align}
N_{[]}(\epsilon, \mathcal F, L_1(\mathbb P_{s,T})) = 2\left(\lfloor 2/\epsilon  \rfloor-1\right) <\infty. \label{eq4-3-1}
\end{align}

For each $j=1,\ldots, 2\lfloor 2/\epsilon  \rfloor -1$, define $\mathcal Q_j =\left\{ q\in \mathcal Q: \tilde y_j \le \|\mathbb P_{s,T} \mathcal T_{q(\cdot)}-P\mathcal T_{q(\cdot)} \|_{\mathcal X} \le \tilde y_{j+1}^- \right\}$. 
Define $q_j^-,q_j^+ \in \mathcal Q_j$ satisfying  
$$\tilde y_j \le \|\mathbb P_{s,T} \mathcal T_{q_j^-(\cdot)}-P\mathcal T_{q_j^-(\cdot)} \|_{\mathcal X}  \le \|\mathbb P_{s,T} \mathcal T_{q(\cdot)}-P\mathcal T_{q(\cdot)} \|_{\mathcal X}  \le \|\mathbb P_{s,T} \mathcal T_{q_j^+(\cdot)}-P\mathcal T_{q_j^+(\cdot)} \|_{\mathcal X}  \le \tilde y_{j+1}^- $$ 
for all $q\in \mathcal Q_j$. 
Then, we have
\begin{align*}
\sup_{\mathcal T_y \in \mathcal F_j} \left\| \mathbb P_{s,T} \mathcal T_y - P \mathcal T_y  \right\|_{\mathcal X} 
&= \sup_{q\in \mathcal Q_j} \left\| \Psi_T(q) - \bar{\Psi}_T(q) - \mathbb E\left[ \Psi_T(q) - \bar{\Psi}_T(q) \right] \right\|_{\mathcal X}\\
&= \left\| \Psi_T(q_j^+) - \bar{\Psi}_T(q_j^+) - \left\{ \Psi_T(q_j^-) - \bar{\Psi}_T(q_j^-) \right\} \right\|_{\mathcal X}\\
&= \delta(q_j^+, q_j^-).
\end{align*}
From Lemma \ref{lem2} and the compactness of $\mathcal Y$, we have $\mathbb E\left[ \delta(q_j^+, q_j^-)^2 \right] =O\left(\frac{s\log s}{T}\right)$, hence, $\left\| \mathbb P_{s,T} - P \right\|_{\mathcal F_j} \stackrel{p}{\to} 0$ is obtained. 
Moreover, by (\ref{eq4-3-1}), 
$$\left\| \mathbb P_{s,T} - P  \right\|_{\mathcal F} = \sup_{q \in \mathcal Q} \left\| \Psi _T (q) - \bar{\Psi}_T(q) - \mathbb E \left[ \Psi _T (q) - \bar{\Psi}_T(q) \right] \right\|_{\mathcal X}\stackrel{p}{\to} 0. $$
From Lemma \ref{lem2}, 
$\sup_{q\in \mathcal Q}\left\| \mathbb E\left[ \Psi_T(q)- \bar{\Psi}_T(q) \right]\right\|_{\mathcal X}\to 0$, as $T\to \infty$, 
which yields $I_1 \stackrel{p}{\to} 0$ from the triangle inequality. 

For $I_2$, we can write by definition
\begin{align*}
\bar{\Psi}_T(q)(x) - \Psi(q)(x) 
%&= \sum_{t=1}^T \alpha_t(x) M_{q(x)}(X_t) - \mathbb E\left[ \psi_{q(x)}(Y_t)|X_t=x \right] \\ 
%& = \sum_{t=1}^T \alpha_t(x) (\tau - F_{\varepsilon}(q(x) - g(X_t))) - (\tau - F_{\varepsilon}(q(x) - g(x)))\\ 
& = \sum_{t=1}^T \alpha_t(x) (F_{\varepsilon}(q(x) - g(x)) - F_{\varepsilon}(q(x) - g(X_t))), 
\end{align*}
for any $q\in \mathcal Q$ and $x \in \mathcal X$. 
From %Assumptions 
\ref{A1}, \ref{A2} %(A-1), (A-2)
 and the Taylor expansion of $F_{\varepsilon}$, there exists some $C>0$ such that
\begin{align*}
\left\|\bar{\Psi}_T - \Psi \right\|_{\mathcal Q} 
%&
\le C \sup_{x\in \mathcal X} \sum_{t=1}^T \alpha_t(x) \left\| x - X_t \right\| %\\
%&
= C \sup_{x\in \mathcal X}  \frac{1}{|\mathcal A_s|} \sum_{A\in \mathcal A_s} \sum_{t\in A^{\mathcal I}} \alpha_{A,t} (x) \left\| x - X_t \right\| %\\
%&\le C \sup_{x\in \mathcal X} \mathrm{diam}(L(x))  \sup_{x'\in \mathcal X}  \frac{1}{|\mathcal A_s|} \sum_{A\in \mathcal A_s} \sum_{t\in A^{\mathcal I}} \alpha_{A,t} (x')\\
%&
= C \left\| \mathrm{diam}(L)\right\|_{\mathcal X}. 
\end{align*}
From this and Corollay 1, we have, $I_2 \stackrel{p}{\to} 0$. \qed

\ \\
\quad By Lemmas \ref{lem1}, \ref{lem2}, and \ref{lem3}, we show Theorem \ref{thm1}.
\paragraph{\textbf{Proof of Theorem \ref{thm1}}}
From (\ref{eq2-2}), 
\begin{align}
\Psi(q_0)(x) = 0  \label{eq4-4-1}
\end{align}
for any $x\in \mathcal X$.
Suppose that a sequence $\{q_T\}\in \mathcal Q$ satisfies $\|\Psi(q_T)\|_{\mathcal X} \to 0$ as $T\to \infty$. 
Then, from (\ref{eq2-1}), we have $\sup_{x\in \mathcal X}|\tau -F_{\epsilon}(q_T(x)-g(x))| \to 0$.  
In addition, from the strictly monotonicity of $F_{\epsilon}$, the inverse function $F_{\epsilon}^{-1}$ is continuous, and since $q_T$ and $g$ are bounded, it follows $\sup_{x\in \mathcal X}|F_{\epsilon}^{-1}(\tau) - q_T(x) +g(x) | \to 0$. 
Notably, $q_0(x) = g(x)+F_{\epsilon}^{-1}(\tau)$, we have
\begin{align}
\|q_T-q_0\|_{\mathcal X} \to 0 \quad \mathrm{as} \ T\to \infty\label{eq4-4-2}
\end{align}
(Identifiability condition). 

From \ref{A7}, there exists $C>0$ such that
\begin{align*}
\|\Psi_T(\hat{q}_T)\|_{\mathcal X} 
= \sup_{x\in \mathcal X} \left| \sum_{t=1}^T \alpha_t(x) \psi_{\hat{q}_T(x)}(Y_t) \right|
%\le C \max_{t\in \{1,\ldots, T\}} \sup_{x\in \mathcal X}  \alpha_t(x)
\le C \max_{t\in \{1,\ldots, T\}} \sup_{x\in \mathcal X}  \max_{A\in \mathcal A_s} \alpha_{A,t}(x),
\end{align*}
and by the definition, 
\begin{align*}
\mathbb E\left[ \alpha_{A,t} (x) \right] \le \frac{1}{k}\mathbb E\left[ \bm 1_{\{X_t \in L(x)\}} \right]=\frac{1}{k|A^{\mathcal I}|}\sum_{t\in A^{\mathcal I}} \mathbb E\left[ \bm 1_{\{X_t \in L(x)\}} \right]=O\left( s^{-1} \right),
\end{align*}
for any $A\in \mathcal A_s$ and $x\in \mathcal X$, which implies 
\begin{align}
\|\Psi_T(\hat{q}_T)\|_{\mathcal X} =O_p\left(s^{-1} \right)=o_p(1). \label{eq4-4-3}
\end{align}

From (\ref{eq4-4-1}), (\ref{eq4-4-2}), (\ref{eq4-4-3})  and 
\begin{align*}
\|\Psi_T- \Psi\|_{\mathcal Q}  
= \sup_{q\in \mathcal Q} \left\| \Psi_T(q) - \Psi(q) \right\|_{\mathcal X}
\stackrel{p}{\to}0
\end{align*}
by Lemma \ref{lem3}, we obtain $\|\hat{q}_T- q_0\|_{\mathcal X} \stackrel{p}{\to} 0$ from Theorem 2.10 of Kosorok (2008)%\cite{kosorok2008}
. \qed

\paragraph{\textbf{Proof of Theorem \ref{thm2}}}
Suppose that for a fixed $s$, $T$ and $B$, $W= (W_{A})_{A\in \mathcal A_s}$ is an multinomial random $|\mathcal A_s|$ vector taking values on $\{0,1,\ldots, B\}^{|\mathcal A_s|}$ with probabilities $1/|\mathcal A_s|$ and number of trials $B$, and which is independent of the data $\mathcal D_T$ and $\xi$. 
Note that $\sum_{A\in \mathcal A_s} W_{A} = B$.   
Then, for any $q\in \mathcal Q$ and $x\in \mathcal X$, we can write  
\begin{align*}
\Psi_T^B(q)(x) 
= \frac{1}{B} \sum_{A\in \mathcal A_s} W_{A} \mathcal T(q;\mathcal I_s, \mathcal J_s,\xi)(x), %\quad
%\end{align*} 
%which implies 
%\begin{align*}
%\alpha_t^B(x) = \frac{1}{B} \sum_{A\in \mathcal A_s} W_A \alpha_{A,t} (x),  
\end{align*}
with 
\begin{align*}
\mathbb E\left[ W_A \right] = \frac{B}{|\mathcal A_s|}, \quad \mathrm {Cov}\left( W_{A_1}, W_{A_2} \right) = \left\{
  \begin{array}{cc}
    \frac{B}{|\mathcal A_s|}\left(1- \frac{1}{|\mathcal A_s|}\right)   & \mathrm{if}\ A_1=A_2   \\
    - \frac{B}{|\mathcal A_s|^2}  & \mathrm{if}\ A_1\neq A_2   \\
  \end{array}
\right.
\end{align*}
Let $\mathbb E_W$ denote taking the expectation over $W= (W_{A})_{A\in \mathcal A_s}$. 
Then, we have
\begin{align*}
\mathbb E_{W}\left[ \Psi_T^B(q)(x) \right] 
&= \Psi_T(q)(x),\\
\mathbb E_{W}\left[ \Psi_T^B(q)(x)^2 \right] 
&= \frac{1}{B|\mathcal A_s|} \sum_{A\in \mathcal A_s}\mathcal T(q;\mathcal I_s, \mathcal J_s,\xi)(x)^2 + \left(1-\frac{1}{B}\right) \Psi_T(q)(x)^2, 
\end{align*} 
which implies that
\begin{align*}
\mathbb E \left[ \left(\Psi_T^B(q)(x) - \Psi_T(q)(x) \right)^2\right] 
&= \frac{1}{B} \left\{ \frac{1}{|\mathcal A_s|} \sum_{A\in \mathcal A_s} \mathbb E\left[ \mathcal T(q;\mathcal I_s, \mathcal J_s,\xi)(x)^2 \right] - \mathbb E \left[ \Psi_T(q)(x) ^2 \right] \right\} = O\left(\frac{1}{B} \right). 
\end{align*} 
By using this, under $B^{-1}=o(1)$, we have $\left\| \Psi_T^B - \Psi_T \right\|_{\mathcal Q} \stackrel{p}{\to} 0$, and from Lemma \ref{lem3},  $\left\| \Psi_T^B - \Psi \right\|_{\mathcal Q} \stackrel{p}{\to} 0$. 
By the same argument of Theorem \ref{thm1}, we obtain $\left\| \hat{q}_T^B - q_0 \right\|_{\mathcal X} \stackrel{p}{\to} 0$  as $T\to \infty$. \qed

\end{document}